\documentclass[11pt,a4paper]{article}

\usepackage{soul}
\usepackage{tikz}
\usepackage{amsthm}
\usepackage{amsmath}
\usepackage{amssymb}
\usepackage{mathrsfs}
\usepackage{geometry}
\usepackage{graphicx}
\usepackage{amsfonts}
\usepackage{epstopdf}
\usepackage{placeins}
\usepackage{enumerate}
\usepackage{enumitem}
\usepackage{hyperref}
\usepackage{subcaption}

\usetikzlibrary{calc}

\numberwithin{equation}{section}

\allowdisplaybreaks

\newcommand{\Rmnum}[1]{\uppercase\expandafter{\romannumeral#1}} % Uppercase roman number

\def\Xint#1{\mathchoice
{\XXint\displaystyle\textstyle{#1}}%
{\XXint\textstyle\scriptstyle{#1}}%
{\XXint\scriptstyle\scriptscriptstyle{#1}}%
{\XXint\scriptscriptstyle\scriptscriptstyle{#1}}%
\!\int}
\def\XXint#1#2#3{{\setbox0=\hbox{$#1{#2#3}{\int}$ }
\vcenter{\hbox{$#2#3$ }}\kern-.6\wd0}}

\def\dashint{\Xint-}

\theoremstyle{plain}
\newtheorem{theorem}{Theorem}[section]
\newtheorem{proposition}[theorem]{Proposition}
\newtheorem{lemma}[theorem]{Lemma}
\newtheorem{corollary}[theorem]{Corollary}

\theoremstyle{definition}
\newtheorem{remark}[theorem]{Remark}

\renewcommand{\thefootnote}{}

\AtEndDocument{
\bigskip{\footnotesize%
   \textsc{Department of Mathematics, Aarhus University, 8000 Aarhus C, Denmark} \par
  \textit{E-mail address}: \texttt{yangmengqh@gmail.com}
}}

\begin{document}

\title{{On singularity of $p$-energy measures on metric measure spaces}}
\author{Meng Yang}
\date{}

\maketitle

\abstract{For $p>1$, we prove that, for a $p$-energy on a volume doubling metric measure space, the Poincar\'e inequality and the cutoff Sobolev inequality, both with $p$-walk dimension strictly larger than $p$, imply that the associated $p$-energy measure is singular with respect to the underlying measure. Under the slow volume regularity condition, we further prove that these two inequalities are equivalent to the resistance estimate; in particular, as part of the proof, we give a simple and direct derivation of the cutoff Sobolev inequality from the Poincar\'e inequality and the capacity upper bound. As a direct corollary, for a large family of fractals and metric measure spaces, including the Sierpi\'nski gasket and the Sierpi\'nski carpet, the $p$-energy measure is singular with respect to the underlying measure for any $p$ strictly greater than the Ahlfors regular conformal dimension.}

\footnote{\textsl{Date}: \today}
\footnote{The author was very grateful to Ryosuke Shimizu for pointing out the recent paper \cite{Sas26} on (canonical) $p$-energy measures and for valuable communications on nonlinear potential theory, and to Aobo Chen for several useful comments.}

\renewcommand{\thefootnote}{\arabic{footnote}}
\setcounter{footnote}{0}

\section{Introduction}

On a large family of fractals, including the Sierpi\'nski gasket and the Sierpi\'nski carpet, there exists a diffusion with a heat kernel satisfying the following two-sided sub-Gaussian estimate
\begin{align*}\label{eq_HKbeta}\tag*{$\text{HK}(\beta)$}
\frac{C_1}{V(x,t^{1/\beta})}\exp\left(-C_2\left(\frac{d(x,y)}{t^{1/\beta}}\right)^{\frac{\beta}{\beta-1}}\right)\le p_t(x,y)\le\frac{C_3}{V(x,t^{1/\beta})}\exp\left(-C_4\left(\frac{d(x,y)}{t^{1/\beta}}\right)^{\frac{\beta}{\beta-1}}\right),\nonumber
\end{align*}
where $\beta$ is a new parameter called the walk dimension, which is always strictly greater than 2 on fractals. For example, $\beta=\log5/\log2$ on the Sierpi\'nski gasket (see \cite{BP88,Kig89}), $\beta\approx2.09697$ on the Sierpi\'nski carpet (see \cite{BB89,BB90,BBS90,BB92,KZ92,HKKZ00}). For $\beta=2$, \ref{eq_HKbeta} is indeed the classical Gaussian estimate, which is equivalent to the conjunction of the volume doubling condition and the Poincar\'e inequality on any complete non-compact Riemannian manifold, see \cite{Gri92,Sal92,Sal95}.

By the standard Dirichlet form theory, a diffusion corresponds to a local regular Dirichlet form associated with an energy measure, see \cite{FOT11} for more details. It was proved in \cite{KM20} that under \ref{eq_HKbeta}, the energy measure is singular with respect to the underlying measure if and only if $\beta>2$, see also \cite{Kus89,Kus93LNM,BST99,Hin05,HN06} for previous results on self-similar sets.

The Dirichlet form framework generalizes the classical Dirichlet integral $\int_{\mathbb{R}^d}|\nabla f(x)|^2\mathrm{d} x$ in $\mathbb{R}^d$. For general $p>1$, extending the classical $p$-energy $\int_{\mathbb{R}^d}|\nabla f(x)|^p\mathrm{d} x$ in $\mathbb{R}^d$, as initiated by \cite{HPS04}, the study of $p$-energy on fractals and general metric measure spaces has been recently advanced considerably, see \cite{CGQ22,Shi24,BC23,MS25,Kig23,CGYZ26,AB25,AES25a}. In this setting, a new parameter $\beta_p$, called the $p$-walk dimension, naturally arises in connection with a $p$-energy. Notably, $\beta_2$ coincides with $\beta$ in \ref{eq_HKbeta}. Analogously, a $p$-energy is also associated with a (canonical) $p$-energy measure, see \cite{Sas26}. It is therefore reasonable to formulate the following conjecture, as was done in \cite[Problem 10.5]{MS25}.

\vspace{5pt}
\noindent{\textbf{Conjecture:}} For a $p$-energy with $p$-walk dimension $\beta_p$, the associated $p$-energy measure is singular with respect to the underlying measure for any $\beta_p>p$.
\vspace{5pt}

The main motivation of this paper is to make progress toward this conjecture. Recall that in the case $p=2$, the main tools used in \cite{KM20} are indeed the Poincar\'e inequality and the cutoff Sobolev inequality both with walk dimension $\beta$, whose conjunction is known to be equivalent to \ref{eq_HKbeta}, see \cite{GHL15}. Recently, a $p$-version of the cutoff Sobolev inequality was introduced by the author in \cite{Yan25a}. Combined with a $p$-version of the Poincar\'e inequality, both involving the $p$-walk dimension $\beta_p$, we will follow an argument similar to that in \cite{KM20} to establish the singularity of the $p$-energy measure with respect to the underlying measure for any $\beta_p>p$.

The problem reduces to proving the Poincar\'e inequality and the cutoff Sobolev inequality on certain fractals or metric measure spaces. The Poincar\'e inequality is relatively straightforward. Indeed, on many such spaces, it arises naturally as a byproduct of the construction of the $p$-energy. The main difficulty lies in the cutoff Sobolev inequality. Even in the case $p=2$, a direct proof of the cutoff Sobolev inequality has long been regarded as highly non-trivial, and its validity is typically established via conditions equivalent to heat kernel estimates.

Recently, in \cite{Yan25b}, the author gave a direct proof of the cutoff Sobolev inequality for $p=2$ on the Sierpi\'nski gasket. In the present paper, we adopt the same approach to show that, under the slow volume regularity condition, which generalizes the strongly recurrent condition from the case $p=2$ to general $p>1$ and is equivalent to $p$ being strictly greater than the Ahlfors regular conformal dimension on homogeneous metric measure spaces, the Poincar\'e inequality and the cutoff Sobolev inequality are equivalent to the resistance estimate. As part of the proof, we give a simple and direct derivation of the cutoff Sobolev inequality from the Poincar\'e inequality and the capacity upper bound. A similar derivation was very recently obtained in \cite{Eri26} under the more general assumption of volume doubling, using a Whitney covering technique. Our proof is direct and may be of independent interest.

Since the resistance estimate has been established on certain homogeneous metric measure spaces, as shown in \cite{Kig23,KS24a}, we can provide a partial answer to the above conjecture that the singularity holds for all $p$ strictly greater than the Ahlfors regular conformal dimension.

Throughout this paper, $p\in(1,+\infty)$ is fixed. The letters $C$, $C_1$, $C_2$, $C_A$, $C_B$ will always refer to some positive constants and may change at each occurrence. The sign $\asymp$ means that the ratio of the two sides is bounded from above and below by positive constants. The sign $\lesssim$ ($\gtrsim$) means that the LHS is bounded by positive constant times the RHS from above (below). We use $\#A$ to denote the cardinality of a set $A$.

\section{Statement of main results}

We say that a function $\Phi:[0,+\infty)\to[0,+\infty)$ is doubling if $\Phi$ is a homeomorphism, which implies that $\Phi$ is strictly increasing continuous and $\Phi(0)=0$, and there exists $C_\Phi>1$, called a doubling constant of $\Phi$, such that $\Phi(2r)\le C_\Phi\Phi(r)$ for any $r>0$. Throughout this paper, we always assume that $\Phi, \Psi$ are two doubling functions with doubling constants $C_\Phi, C_\Psi$, respectively.

Let $(X,d,m)$ be a \emph{complete} metric measure space, that is, $(X,d)$ is a complete locally compact separable metric space and $m$ is a positive Radon measure on $X$ with full support. Throughout this paper, we always assume that all metric balls are relatively compact. For any $x\in X$ and any $r>0$, denote $B(x,r)=\{y\in X:d(x,y)<r\}$ and $V(x,r)=m(B(x,r))$. If $B=B(x,r)$, then denote $\delta B=B(x,\delta r)$ for any $\delta>0$. Let $\mathrm{diam}(X)=\sup\{d(x,y):x,y\in X\}$ be the diameter of $(X,d)$. Let $\mathcal{B}(X)$ be the family of all Borel measurable subsets of $X$. Let $C(X)$ be the family of all continuous functions on $X$. Let $C_c(X)$ be the family of all continuous functions on $X$ with compact support. Denote $\dashint_A=\frac{1}{m(A)}\int_A$ and $u_A=\dashint_Au\mathrm{d} m$ for any measurable set $A$ with $m(A)\in(0,+\infty)$ and any function $u$ such that the integral $\int_Au\mathrm{d} m$ is well-defined.

Let $\varepsilon>0$. We say that $W$ is an $\varepsilon$-net (of $(X,d)$) if $W\subseteq X$ satisfies that for any distinct $x,y\in W$, we have $d(x,y)\ge \varepsilon$, and for any $z\in X$, there exists $x\in W$ such that $d(x,z)< \varepsilon$. Since $(X,d)$ is separable, all $\varepsilon$-nets are countable.

We say that the chain condition \ref{eq_CC} holds if there exists $C_{cc}>0$ such that for any $x,y\in X$, for any positive integer $n$, there exists a sequence $\{x_k:0\le k\le n\}$ of points in $X$ with $x_0=x$ and $x_n=y$ such that
\begin{equation*}\label{eq_CC}\tag*{$\text{CC}$}
d(x_k,x_{k-1})\le C_{cc} \frac{d(x,y)}{n}\text{ for any }k=1,\ldots,n.
\end{equation*}
Throughout this paper, we always assume \ref{eq_CC}.

We say that the volume doubling condition \ref{eq_VD} holds if there exists $C_{VD}>0$ such that
\begin{equation*}\label{eq_VD}\tag*{$\text{VD}$}
V(x,2r)\le C_{VD}V(x,r)\text{ for any }x\in X,r>0.
\end{equation*}

We say that the volume regular condition \ref{eq_VPhi} holds if there exists $C_{VR}>0$ such that
\begin{equation*}\label{eq_VPhi}\tag*{$\text{V}(\Phi)$}
\frac{1}{C_{VR}}\Phi(r)\le V(x,r)\le C_{VR}\Phi(r)\text{ for any }x\in X,r\in(0,\mathrm{diam}(X)).
\end{equation*}
For $d_h>0$, we say that the Ahlfors regular condition \hypertarget{eq_Vdh}{$\text{V}(d_h)$} holds if \ref{eq_VPhi} holds with $\Phi:r\mapsto r^{d_h}$.

We say that the slow volume regular condition \ref{eq_SVR} holds if \ref{eq_VPhi} holds and there exists an increasing continuous function $\Theta:[0,1]\to[0,+\infty)$ with $\Theta(0)=0$, $\Theta(1)\ge1$ and $\sum_{n=0}^{+\infty}\Theta(\frac{1}{2^n})^{1/p}<+\infty$ such that
\begin{equation*}\label{eq_SVR}\tag*{$\text{SVR}(\Phi,\Psi)$}
\frac{\Phi(R)}{\Phi(r)}\le\Theta\left(\frac{r}{R}\right)\frac{\Psi(R)}{\Psi(r)}
\text{ for any }R,r\in(0,\mathrm{diam}(X))\text{ with }r\le R.
\end{equation*}
For $d_h,\beta_p>0$, we say that the slow volume regular condition \hypertarget{eq_SVRdb}{$\text{SVR}(d_h,\beta_p)$} holds if \hyperlink{eq_Vdh}{$\text{V}(d_h)$} holds and $d_h<\beta_p$.

We say that $(\mathcal{E},\mathcal{F})$ is a $p$-energy on $(X,d,m)$ if $\mathcal{F}$ is a dense subspace of $L^p(X;m)$ and $\mathcal{E}:\mathcal{F}\to[0,+\infty) $ satisfies the following conditions.

\begin{enumerate}[label=(\arabic*)]
\item $\mathcal{E}^{1/p}$ is a semi-norm on $\mathcal{F}$, that is, for any $f,g\in\mathcal{F}$, $c\in\mathbb{R}$, we have $\mathcal{E}(f)\ge0$, $\mathcal{E}(cf)^{1/p}=|c|\mathcal{E}(f)^{1/p}$ and $\mathcal{E}(f+g)^{1/p}\le\mathcal{E}(f)^{1/p}+\mathcal{E}(g)^{1/p}$.
\item (Closed property) $(\mathcal{F},\mathcal{E}(\cdot)^{1/p}+\lVert {\cdot}\rVert_{L^p(X;m)})$ is a Banach space.
\item (Markovian property) For any $\varphi\in C(\mathbb{R})$ with $\varphi(0)=0$ and $|\varphi(t)-\varphi(s)|\le|t-s|$ for any $t$, $s\in\mathbb{R}$, for any $f\in\mathcal{F}$, we have $\varphi(f)\in\mathcal{F}$ and $\mathcal{E}(\varphi(f))\le\mathcal{E}(f)$.
\item (Regular property) $\mathcal{F}\cap C_c(X)$ is uniformly dense in $C_c(X)$ and $(\mathcal{E}(\cdot)^{1/p}+\lVert {\cdot}\rVert_{L^p(X;m)})$-dense in $\mathcal{F}$.
\item (Strongly local property) For any $f,g\in\mathcal{F}$ with compact support and $g$ constant in an open neighborhood of $\mathrm{supp}(f)$, we have $\mathcal{E}(f+g)=\mathcal{E}(f)+\mathcal{E}(g)$.
\item ($p$-Clarkson's inequality) For any $f,g\in\mathcal{F}$, we have
\begin{equation*}\label{eq_Cla}\tag*{$\text{Cla}$}
\begin{cases}
\mathcal{E}(f+g)+\mathcal{E}(f-g)\ge2 \left(\mathcal{E}(f)^{\frac{1}{p-1}}+\mathcal{E}(g)^{\frac{1}{p-1}}\right)^{p-1}&\text{if }p\in(1,2],\\
\mathcal{E}(f+g)+\mathcal{E}(f-g)\le2 \left(\mathcal{E}(f)^{\frac{1}{p-1}}+\mathcal{E}(g)^{\frac{1}{p-1}}\right)^{p-1}&\text{if }p\in[2,+\infty).\\
\end{cases}
\end{equation*}
\end{enumerate}
Moreover, we also always assume the following condition.
\begin{itemize}
\item ($\mathcal{F}\cap L^\infty(X;m)$ is an algebra) For any $f,g\in\mathcal{F}\cap L^\infty(X;m)$, we have $fg\in\mathcal{F}$ and
\begin{equation*}\label{eq_Alg}\tag*{$\text{Alg}$}
\mathcal{E}(fg)^{{1}/{p}}\le \lVert {f}\rVert_{L^\infty(X;m)}\mathcal{E}(g)^{1/p}+\lVert {g}\rVert_{L^\infty(X;m)}\mathcal{E}(f)^{1/p}.
\end{equation*}
\end{itemize}
Denote $\mathcal{E}_1(\cdot)=\mathcal{E}(\cdot)+\lVert {\cdot}\rVert^p_{L^p(X;m)}$. Indeed, a general condition called the generalized $p$-contraction property was introduced in \cite{KS24a}, which implies \ref{eq_Cla}, \ref{eq_Alg}, and holds on a large family of metric measure spaces.

We list some basic properties of $p$-energies as follows. For any $f,g\in\mathcal{F}$, the derivative
$$\mathcal{E}(f;g)=\frac{1}{p}\frac{\mathrm{d}}{\mathrm{d} t}\mathcal{E}(f+tg)|_{t=0}\in\mathbb{R}$$
exists, the map $\mathcal{E}(f;\cdot):\mathcal{F}\to\mathbb{R}$ is linear, $\mathcal{E}(f;f)=\mathcal{E}(f)$. Moreover, for any $f$, $g\in\mathcal{F}$, for any $a\in\mathbb{R}$, we have
\begin{equation}\label{eq_quasi_strict}
\mathbb{R}\ni t\mapsto\mathcal{E}(f+tg;g)\in\mathbb{R}\text{ is strictly increasing if and only if }\mathcal{E}(g)>0,
\end{equation}
$$\mathcal{E}(af;g)=\mathrm{sgn}(a)|a|^{p-1}\mathcal{E}(f;g),$$
$$|\mathcal{E}(f;g)|\le\mathcal{E}(f)^{(p-1)/p}\mathcal{E}(g)^{1/p}.$$
Moreover, all of the above results remain valid with $\mathcal{E}$ replaced by $\mathcal{E}_1$, and for any $f,g\in\mathcal{F}$, we have
$$\mathcal{E}_1(f;g)=\mathcal{E}(f;g)+\int_X\mathrm{sgn}(f)|f|^{p-1}g\mathrm{d} m.$$
See \cite[Theorem 3.7, Corollary 3.25]{KS24a} for the proofs of these results.

By \cite[Theorem 2.4]{Sas26}, a $p$-energy $(\mathcal{E},\mathcal{F})$ corresponds to a (canonical) $p$-energy measure $\Gamma:\mathcal{F}\times\mathcal{B}(X)\to[0,+\infty)$, $(f,A)\mapsto\Gamma(f)(A)$ satisfying the following conditions.
\begin{enumerate}[label=(\alph*),ref=(\alph*)]
\item\label{item_meas1} For any $f\in\mathcal{F}$, $\Gamma(f)(\cdot)$ is a positive Radon measure on $X$ with $\Gamma(f)(X)=\mathcal{E}(f)$.
\item\label{item_meas2} For any $A\in\mathcal{B}(X)$, $\Gamma(\cdot)(A)^{1/p}$ is a semi-norm on $\mathcal{F}$.
\item\label{item_meas3} For any $f,g\in\mathcal{F}\cap C_c(X)$, $A\in\mathcal{B}(X)$, if $f-g$ is constant on $A$, then $\Gamma(f)(A)=\Gamma(g)(A)$.
\item\label{item_meas4} ($p$-Clarkson's inequality) For any $f,g\in\mathcal{F}$, for any $A\in\mathcal{B}(X)$, we have
\begin{equation*}
\begin{cases}
\Gamma(f+g)(A)+\Gamma(f-g)(A)\ge2 \left(\Gamma(f)(A)^{\frac{1}{p-1}}+\Gamma(g)(A)^{\frac{1}{p-1}}\right)^{p-1}&\text{if }p\in(1,2],\\
\Gamma(f+g)(A)+\Gamma(f-g)(A)\le2 \left(\Gamma(f)(A)^{\frac{1}{p-1}}+\Gamma(g)(A)^{\frac{1}{p-1}}\right)^{p-1}&\text{if }p\in[2,+\infty).\\
\end{cases}
\end{equation*}
\item\label{item_meas5} (Chain rule) For any $f,g\in\mathcal{F}\cap C_c(X)$, for any piecewise $C^1$ functions $\varphi,\psi:\mathbb{R}\to\mathbb{R}$ with $\varphi(0)=\psi(0)=0$, we have
$$\mathrm{d}\Gamma\left(\varphi(f);\psi(g)\right)=\mathrm{sgn}(\varphi'(f))|\varphi'(f)|^{p-1}\psi'(g)\mathrm{d}\Gamma(f;g),$$
where $\Gamma(f;g)$ is a signed measure given by $\Gamma(f;g)=\frac{1}{p}\frac{\mathrm{d}}{\mathrm{d}t}\Gamma(f+tg)|_{t=0}$.
\item\label{item_meas6} (Leibniz rule) For any $f,g,h\in \mathcal{F}\cap C_c(X)$, we have $\mathrm{d}\Gamma(f;gh)=g \mathrm{d}\Gamma(f;h)+h \mathrm{d}\Gamma(f;g)$.
\end{enumerate}
Using the chain rule, we have the following condition.
\begin{itemize}
\item (Strong sub-additivity) For any $f,g\in\mathcal{F}$, we have $f\vee g, f\wedge g\in\mathcal{F}$ and
\begin{equation*}\label{eq_SubAdd}\tag*{$\text{SubAdd}$}
\mathcal{E}(f\vee g)+\mathcal{E}(f\wedge g)\le\mathcal{E}(f)+\mathcal{E}(g).
\end{equation*}
\end{itemize}

Let
\begin{align*}
&\mathcal{F}_{\mathrm{loc}}=\left\{u:
\begin{array}{l}
\text{for any relatively compact open set }U,\\
\text{there exists }u^\#\in\mathcal{F}\text{ such that }u=u^\#\text{ }m\text{-a.e. in }U
\end{array}
\right\}.
\end{align*}
For any $u\in \mathcal{F}_{\mathrm{loc}}$, let $\Gamma(u)|_U=\Gamma(u^\#)|_U$, where $u^\#, U$ are given as above, then $\Gamma(u)$ is a well-defined positive Radon measure on $X$, as follows from \ref{item_meas3} and \ref{item_meas6} together with an argument similar to that in \cite[Corollary 3.2.1]{FOT11}.

We say that the Poincar\'e inequality \ref{eq_PI} holds if there exist $C_{PI}>0$, $A_{PI}\ge1$ such that for any ball $B$ with radius $r\in(0,\mathrm{diam}(X)/A_{PI})$, for any $f\in\mathcal{F}$, we have
\begin{equation*}\label{eq_PI}\tag*{$\text{PI}(\Psi)$}
\int_B\lvert f-f_B\rvert^p\mathrm{d} m\le C_{PI}\Psi(r)\int_{A_{PI}B}\mathrm{d}\Gamma(f).
\end{equation*}
For $\beta_p>0$, we say that the Poincar\'e inequality \hypertarget{eq_PIbeta}{$\text{PI}(\beta_p)$} holds if \ref{eq_PI} holds with $\Psi:r\mapsto r^{\beta_p}$.

Let $U$, $V$ be two open subsets of $X$ satisfying $U\subseteq\overline{U}\subseteq V$. We say that $\phi\in\mathcal{F}$ is a cutoff function for $U\subseteq V$ if $0\le\phi\le1$ in $X$, $\phi=1$ in an open neighborhood of $\overline{U}$ and $\mathrm{supp}(\phi)\subseteq V$, where $\mathrm{supp}(f)$ refers to the support of the measure of $|f|\mathrm{d} m$ for any given function $f$.

We say that the cutoff Sobolev inequality \ref{eq_CS} holds if there exist $C_{1}, C_{2}>0$, $A_{S}>1$ such that for any ball $B(x,r)$, there exists a cutoff function $\phi\in\mathcal{F}$ for $B(x,r)\subseteq B(x,A_Sr)$ such that for any $f\in\mathcal{F}$, we have
\begin{equation*}\label{eq_CS}\tag*{$\text{CS}(\Psi)$}
\int_{B(x,A_{S}r)}|\widetilde{f}|^p \mathrm{d}\Gamma(\phi)\le C_{1}\int_{B(x,A_{S}r)}\mathrm{d}\Gamma(f)+\frac{C_{2}}{\Psi(r)}\int_{B(x,A_{S}r)}|f|^p\mathrm{d} m,
\end{equation*}
where $\widetilde{f}$ is a quasi-continuous modification of $f$, such that $\widetilde{f}$ is uniquely determined $\Gamma(\phi)$-a.e. in $X$, see \cite[Section 8]{Yan25a} for more details. For $\beta_p>0$, we say that the cutoff Sobolev inequality \hypertarget{eq_CSbeta}{$\text{CS}(\beta_p)$} holds if \ref{eq_CS} holds with $\Psi:r\mapsto r^{\beta_p}$.

Our first main result is as follows.

\begin{theorem}\label{thm_singular}
Assume \ref{eq_VD}, \ref{eq_PI}, \ref{eq_CS} and
\begin{equation}\label{eq_Psi_singular}
\varliminf_{\lambda\to+\infty}\varliminf_{r\downarrow0}\frac{\lambda^p\Psi\left(\frac{r}{\lambda}\right)}{\Psi(r)}=0,
\end{equation}
then for any $f\in\mathcal{F}$, we have $\Gamma(f)\perp m$. In particular, assume \ref{eq_VD}, \hyperlink{eq_PIbeta}{$\text{PI}(\beta_p)$}, \hyperlink{eq_CSbeta}{$\text{CS}(\beta_p)$} with $\beta_p>p$, then for any $f\in\mathcal{F}$, we have $\Gamma(f)\perp m$.
\end{theorem}

\begin{remark}
We will follow an argument from \cite{KM20}, where the case $p=2$ was considered.
\end{remark}

Let $A_1,A_2\in\mathcal{B}(X)$. We define the capacity between $A_1$, $A_2$ as
\begin{align*}
&\mathrm{cap}(A_1,A_2)=\inf\left\{\mathcal{E}(\varphi):\varphi\in\mathcal{F},
\begin{array}{l}
\varphi=1\text{ in an open neighborhood of }A_1,\\
\varphi=0\text{ in an open neighborhood of }A_2
\end{array}
\right\},
\end{align*}
here we use the convention that $\inf\emptyset=+\infty$.

We say that the capacity upper bound \ref{eq_cap} holds if there exist $C_{cap}>0$, $A_{cap}>1$ such that for any ball $B(x,r)$, we have
\begin{equation*}\label{eq_cap}\tag*{$\text{cap}(\Psi)_{\le}$}
\mathrm{cap}\left(B(x,r),X\backslash B(x,A_{cap}r)\right)\le C_{cap} \frac{V(x,r)}{\Psi(r)}.
\end{equation*}
For $\beta_p>0$, we say that the capacity upper bound \hypertarget{eq_capbeta}{$\text{cap}(\beta_p)_{\le}$} holds if \ref{eq_cap} holds with $\Psi:r\mapsto r^{\beta_p}$. Under \ref{eq_VD}, by taking $f\equiv1$ in $B(x,A_Sr)$, it is easy to see that \ref{eq_CS} (resp. \hyperlink{eq_CSbeta}{$\text{CS}(\beta_p)$}) implies \ref{eq_cap} (resp. \hyperlink{eq_capbeta}{$\text{cap}(\beta_p)_{\le}$}).

Let $A,B\in\mathcal{B}(X)$. The effective resistance between $A$ and $B$ is defined as follows.
$$R(A,B)=\left(\inf\left\{\mathcal{E}(u):u\in\mathcal{F},u=0\text{ on }A,u=1\text{ on }B\right\}\right)^{-1},$$
here we use the convention that $0^{-1}=+\infty$  and $(+\infty)^{-1}=0$. By definition, we have
$$R(A,B)=\mathrm{cap}(A,B)^{-1}.$$
For $x,y\in X$, write
$$R(x,A)=R(\{x\},A),R(x,y)=R(\{x\},\{y\}),$$
then
$$R(x,y)=\sup\left\{\frac{|u(x)-u(y)|^p}{\mathcal{E}(u)}:u\in\mathcal{F},\mathcal{E}(u)>0\right\},$$
which implies that
$$|u(x)-u(y)|^p\le R(x,y)\mathcal{E}(u)\text{ for any }u\in\mathcal{F},x,y\in X.$$
If $A_1\subseteq A_2$, $B_1\subseteq B_2$, then it is obvious that
$$R(A_1,B_1)\ge R(A_2,B_2).$$

We say that the resistance estimate \ref{eq_res} holds if there exists $C_R>0$ such that for any $x,y\in X$ with $x\ne y$, we have
\begin{equation*}\label{eq_res}\tag*{$\text{R}(\Phi,\Psi)$}
\frac{1}{C_R}\frac{\Psi(d(x,y))}{\Phi(d(x,y))}\le R(x,y)\le C_R\frac{\Psi(d(x,y))}{\Phi(d(x,y))}.
\end{equation*}
For $d_h, \beta_p>0$, we say that the resistance estimate \hypertarget{eq_Rdb}{$\text{R}(d_h,\beta_p)$} holds if \ref{eq_res} holds with $\Phi:r\mapsto r^{d_h}$, $\Psi:r\mapsto r^{\beta_p}$.

Our second main result is as follows.

\begin{theorem}\label{thm_equiv}
Assume \ref{eq_SVR}. The followings are equivalent.
\begin{enumerate}[label=(\arabic*)]
\item \ref{eq_res}.
\item \ref{eq_PI} and \ref{eq_cap}.
\item \ref{eq_PI} and \ref{eq_CS}.
\end{enumerate}
\end{theorem}

\begin{remark}
The key novelty of this result lies in the fact that, under \ref{eq_SVR}, the conjunction of \ref{eq_PI} and \ref{eq_cap} implies \ref{eq_CS}. We employ an idea recently introduced by the author in \cite{Yan25b}.
\end{remark}

As an application of our results, we consider $p$-energies on some concrete fractals. On the Sierpi\'nski gasket, let $d_h=\frac{\log3}{\log2}$ be the Hausdorff dimension and $m$ the $d_h$-dimensional Hausdorff measure, let $\mathrm{dim}_{\mathrm{ARC}}=1$ be the Ahlfors regular conformal dimension, see \cite{TW06,Car14}. By \cite{HPS04,CGQ22}, for any $p>1$, there exists a $p$-energy $(\mathcal{E},\mathcal{F})$ with $p$-walk dimension $\beta_p>d_h$. On the Sierpi\'nski carpet, let $d_h=\frac{\log8}{\log3}$ be the Hausdorff dimension and $m$ the $d_h$-dimensional Hausdorff measure, let $\mathrm{dim}_{\mathrm{ARC}}$ be the Ahlfors regular conformal dimension which lies in $[1+\frac{\log2}{\log3},\frac{\log8}{\log3})$, see \cite[Example 4.3.1]{MT10book}, or \cite{Tys00a} and \cite{KL04}. By \cite{Shi24,MS25}, for any $p>1$, there exists a $p$-energy $(\mathcal{E},\mathcal{F})$ with $p$-walk dimension $\beta_p$, and $\beta_p>d_h$ if and only if $p>\mathrm{dim}_{\mathrm{ARC}}$.

\begin{corollary}\label{cor_SGSC}
On the Sierpi\'nski gasket and the Sierpi\'nski carpet, for any $p>\mathrm{dim}_{\mathrm{ARC}}$, we have \hyperlink{eq_SVRdb}{$\text{SVR}(d_h,\beta_p)$}, \hyperlink{eq_Rdb}{$\text{R}(d_h,\beta_p)$} hold, hence \hyperlink{eq_PIbeta}{$\text{PI}(\beta_p)$}, \hyperlink{eq_capbeta}{$\text{cap}(\beta_p)_\le$}, \hyperlink{eq_CSbeta}{$\text{CS}(\beta_p)$} all hold, which gives that for any $f\in\mathcal{F}$, we have $\Gamma(f)\perp m$.
\end{corollary}

\begin{proof}
For any $p>\mathrm{dim}_{\mathrm{ARC}}$, since $\beta_p>d_h$, we have \hyperlink{eq_SVRdb}{$\text{SVR}(d_h,\beta_p)$} holds. By \cite[Theorem 4.6]{Kig23} and \cite[Theorem 8.30]{KS24a}, we have \hyperlink{eq_Rdb}{$\text{R}(d_h,\beta_p)$} holds on the Sierpi\'nski gasket. By \cite[Theorem 3.21]{Kig23} and \cite[Theorem 8.19, Remark 8.20]{KS24a}, we have \hyperlink{eq_Rdb}{$\text{R}(d_h,\beta_p)$} holds on the Sierpi\'nski carpet. By Theorem \ref{thm_equiv}, we have \hyperlink{eq_PIbeta}{$\text{PI}(\beta_p)$}, \hyperlink{eq_capbeta}{$\text{cap}(\beta_p)_\le$}, \hyperlink{eq_CSbeta}{$\text{CS}(\beta_p)$} all hold. By \cite[Theorem 9.13, Theorem 9.8]{KS24a}, we have $\beta_p>p$. By Theorem \ref{thm_singular}, we have the singularity result.
\end{proof}

\begin{remark}\hspace{0em}
\begin{enumerate}[label=(\arabic*)]
\item In addition to the Sierpi\'nski gasket, the above results also apply to strongly symmetric p.c.f. self-similar sets, as defined in \cite[Definition 3.8.3, Definition 3.8.4]{Kig01book}. More precisely, under an additional assumption \cite[Assumption B.6]{KS24a}, by \cite[Theorem B.8]{KS24a}, these sets have Ahlfors regular conformal dimension equal to 1. Hence for any $p>1$, we have $\beta_p>d_h$, then by \cite[Theorem 4.6]{Kig23} and \cite[Theorem 8.30]{KS24a}, the above results also hold.
\item In addition to the Sierpi\'nski carpet, the above results also apply to $p$-conductively homogeneous compact metric spaces endowed with $p$-energies, as defined in \cite{Kig20LNM,Kig23}. More precisely, for any $p>\mathrm{dim}_{\mathrm{ARC}}$, we have $\beta_p>d_h$, then by \cite[Theorem 3.21]{Kig23} and \cite[Theorem 8.19]{KS24a}, the above results also hold.
\item Due to the heavy details involved in introducing the above spaces, we restrict our attention to the Sierpi\'nski gasket and the Sierpi\'nski carpet as representative examples to illustrate the applications of our results.
\end{enumerate}
\end{remark}

This paper is organized as follows. In Section \ref{sec_CS}, we give some results related to the cutoff Sobolev inequality. In Section \ref{sec_singular}, we prove Theorem \ref{thm_singular}. In Section \ref{sec_equiv}, we prove Theorem \ref{thm_equiv}. In Section \ref{sec_potential}, we give some necessary results from potential theory.

\section{Cutoff Sobolev inequalities}\label{sec_CS}

In this section, we give the self-improvement property of the cutoff Sobolev inequality, along with a reverse Poincar\'e inequality that follows as a consequence.

We have the self-improvement property of the cutoff Sobolev inequality as follows.

\begin{proposition}\label{prop_CS_self}
Assume \ref{eq_VD}, \ref{eq_CS}. Then for any $\delta>0$, there exists $C_\delta>0$ depending on $\delta$, such that for any $x\in X$, for any $R, r>0$, there exists a cutoff function $\phi\in\mathcal{F}$ for $B(x,R)\subseteq B(x,R+r)$ such that for any $f\in\mathcal{F}$, we have
\begin{align}
&\int_{B(x,R+r)\backslash\overline{B(x,R)}}|\widetilde{f}|^p \mathrm{d}\Gamma(\phi)\nonumber\\
&\le\delta\int_{B(x,R+r)\backslash\overline{B(x,R)}}|\widetilde{\phi}|^p\mathrm{d}\Gamma(f)+\frac{C_\delta}{\Psi(r)}\int_{B(x,R+r)\backslash\overline{B(x,R)}}|f|^p\mathrm{d} m,\label{eq_CS_self}
\end{align}
where $\widetilde{f}, \widetilde{\phi}$ are quasi-continuous modifications of $f, \phi$, respectively.
\end{proposition}

\begin{proof}
We follow the same argument as in \cite[Lemma 5.1]{AB15}, \cite[PROPOSITION 5.11]{BM18} and \cite[Lemma 6.2]{Mur24a}, adapted to the $p$-energy setting. Firstly, we prove that there exists $C>0$ such that for any $x\in X$, for any $R, r>0$, there exists a cutoff function $\phi\in\mathcal{F}$ for $B(x,R)\subseteq B(x,R+r)$ such that for any $f\in\mathcal{F}$, we have
\begin{align}
&\int_{B(x,R+r)\backslash\overline{B(x,R)}}|\widetilde{f}|^p\mathrm{d}\Gamma(\phi)\nonumber\\
&\le C\int_{B(x,R+r)\backslash\overline{B(x,R)}}\mathrm{d}\Gamma(f)+\frac{C}{\Psi(r)}\int_{B(x,R+r)\backslash\overline{B(x,R)}}|f|^p\mathrm{d} m.\label{eq_CS_self1}
\end{align}
Let $C_1$, $C_2$, $A=A_S$ be the constants in \ref{eq_CS}. Let $L=4A+4$. Let $V$ be an $\frac{r}{L}$-net. For any $v\in V$, by \ref{eq_CS}, there exists a cutoff function $\phi_v\in\mathcal{F}$ for $B(v,\frac{r}{L})\subseteq B(v,\frac{Ar}{L})$ such that
$$\int_{B(v,\frac{Ar}{L})}|\widetilde{f}|^p\mathrm{d}\Gamma(\phi_v)\le C_{1}\int_{B(v,\frac{Ar}{L})}\mathrm{d}\Gamma(f)+\frac{C_{2}}{\Psi(\frac{r}{L})}\int_{B(v,\frac{Ar}{L})}|f|^p\mathrm{d} m.$$
Let ${V}_1=V\cap B(x,R+\frac{r}{2})$, then by \ref{eq_VD}, we have $\#V_1<+\infty$. Let
$$\phi=\max_{v\in{V}_1}\phi_v,$$
then by \ref{eq_SubAdd}, we have $\phi\in\mathcal{F}$. Moreover, $0\le\phi\le1$ in $X$, $\phi=1$ in $B(x,R+\frac{r}{2}-\frac{r}{L})\supseteq\overline{B(x,R)}$, $\mathrm{supp}(\phi)\subseteq B(x,R+\frac{r}{2}+\frac{Ar}{L})\subseteq B(x,R+r)$, hence $\phi\in\mathcal{F}$ is a cutoff function for $B(x,R)\subseteq B(x,R+r)$ and
$$\mathrm{supp}(\Gamma(\phi))\subseteq B(x,R+\frac{r}{2}+\frac{Ar}{L})\backslash{B(x,R+\frac{r}{2}-\frac{r}{L})}.$$
Let $W={V}_1\backslash B(x,R+\frac{r}{2}-\frac{r}{L}-\frac{Ar}{L})$, then
$$\phi=\max_{v\in W}\phi_v\text{ in }B(x,R+\frac{r}{2}+\frac{Ar}{L})\backslash{B(x,R+\frac{r}{2}-\frac{r}{L})},$$
hence by \ref{eq_SubAdd},
\begin{align*}
&\int_{B(x,R+r)\backslash\overline{B(x,R)}}|\widetilde{f}|^p\mathrm{d}\Gamma(\phi)=\int_{B(x,R+\frac{r}{2}+\frac{Ar}{L})\backslash{B(x,R+\frac{r}{2}-\frac{r}{L})}}|\widetilde{f}|^p\mathrm{d}\Gamma(\phi)\\
&\le\sum_{v\in W}\int_{{B(x,R+\frac{r}{2}+\frac{Ar}{L})\backslash{B(x,R+\frac{r}{2}-\frac{r}{L})}}}|\widetilde{f}|^p\mathrm{d}\Gamma(\phi_v)\le\sum_{v\in W}\int_{{B(v,\frac{Ar}{L})}}|\widetilde{f}|^p\mathrm{d}\Gamma(\phi_v)\\
&\le\sum_{v\in W}\left(C_{1}\int_{B(v,\frac{Ar}{L})}\mathrm{d}\Gamma(f)+\frac{C_{2}}{\Psi(\frac{r}{L})}\int_{B(v,\frac{Ar}{L})}|f|^p\mathrm{d} m\right)\\
&=C_1\int_X \left(\sum_{v\in W}1_{B(v,\frac{Ar}{L})}\right)\mathrm{d}\Gamma(f)+\frac{C_{2}}{\Psi(\frac{r}{L})}\int_X\left(\sum_{v\in W}1_{B(v,\frac{Ar}{L})}\right)|f|^p\mathrm{d} m.
\end{align*}
By \ref{eq_VD}, there exists some positive integer $N$ depending only on $C_{VD}, A, L$ such that
$$\sum_{v\in W}1_{B(v,\frac{Ar}{L})}\le N1_{\cup_{v\in W}B(v,\frac{Ar}{L})},$$
where
$${\bigcup_{v\in W}B(v,\frac{Ar}{L})}\subseteq B(x,R+\frac{r}{2}+\frac{Ar}{L})\backslash\overline{B(x,R+\frac{r}{2}-\frac{r}{L}-\frac{Ar}{L}-\frac{Ar}{L})}\subseteq B(x,R+r)\backslash\overline{B(x,R)},$$
hence
\begin{align*}
&\int_{B(x,R+r)\backslash\overline{B(x,R)}}|\widetilde{f}|^p\mathrm{d}\Gamma(\phi)\\
&\le C_1N\int_{B(x,R+r)\backslash\overline{B(x,R)}}\mathrm{d}\Gamma(f)+\frac{C_2N}{\Psi(\frac{r}{L})}\int_{B(x,R+r)\backslash\overline{B(x,R)}}|f|^p\mathrm{d} m\\
&\le C_1N\int_{B(x,R+r)\backslash\overline{B(x,R)}}\mathrm{d}\Gamma(f)+\frac{C_2NC_3}{\Psi({r})}\int_{B(x,R+r)\backslash\overline{B(x,R)}}|f|^p\mathrm{d} m,
\end{align*}
where $C_3$ is some positive constant depending only on $C_\Psi, L$. Therefore, we have (\ref{eq_CS_self1}) with $C=\max\{C_1N,C_2NC_3\}$.

Secondly, we prove (\ref{eq_CS_self}). Let $a\in(0,1)$ be chosen later. Let $s_n=\frac{1-a}{a}a^nr$ for any $n\ge1$, then $\sum_{n=1}^{+\infty}s_n=r$. Let $r_0=0$ and $r_n=\sum_{k=1}^ns_k$ for any $n\ge1$, then $r_n\uparrow r$. For any $n\ge0$, let $B_n=B(x,R+r_n)$, then $B(x,R)=B_0\subseteq B_n\uparrow B(x,R+r)$. For any $n\ge0$, by (\ref{eq_CS_self1}), there exists a cutoff function $\phi_n\in\mathcal{F}$ for $B_n\subseteq B_{n+1}$ such that
$$\int_{B_{n+1}\backslash\overline{B_n}}|\widetilde{f}|^p\mathrm{d}\Gamma(\phi_n)\le C\int_{{B_{n+1}\backslash\overline{B_n}}}\mathrm{d}\Gamma(f)+\frac{C}{\Psi(s_{n+1})}\int_{{B_{n+1}\backslash\overline{B_n}}}|f|^p\mathrm{d} m.$$
Let $b\in(0,1)$ be chosen later. Let
$$\phi=\sum_{k=0}^{+\infty}(b^k-b^{k+1})\phi_k,$$
then $0\le\phi\le1$ in $X$, $\phi=1$ in $B(x,R)$, $\phi=0$ on $X\backslash B(x,R+r)$. For any $n\ge0$, in $B_{n+1}$, we have $\phi=\sum_{k=0}^n(b^k-b^{k+1})\phi_k+b^{n+1}$, hence $\phi\in\mathcal{F}_{\mathrm{loc}}$, in $B_{n+1}\backslash\overline{B_n}$, we have $\phi=(b^{n}-b^{n+1})\phi_{n}+b^{n+1}=(1-b)b^{n}\phi_{n}+b^{n+1}$, hence $b^{n+1}\le\phi\le b^n$ in $B_{n+1}\backslash\overline{B_n}$. Therefore, we have
\begin{align*}
&\int_{B(x,R+r)\backslash\overline{B(x,R)}}|\widetilde{f}|^p\mathrm{d}\Gamma(\phi)=\sum_{n=0}^{+\infty}\int_{B_{n+1}\backslash\overline{B_n}}|\widetilde{f}|^p\mathrm{d}\Gamma(\phi)=\sum_{n=0}^{+\infty}(1-b)^pb^{pn}\int_{B_{n+1}\backslash\overline{B_n}}|\widetilde{f}|^p\mathrm{d}\Gamma(\phi_n)\\
&\le\sum_{n=0}^{+\infty}(1-b)^pb^{pn}\left(C\int_{{B_{n+1}\backslash\overline{B_n}}}\mathrm{d}\Gamma(f)+\frac{C}{\Psi(s_{n+1})}\int_{{B_{n+1}\backslash\overline{B_n}}}|f|^p\mathrm{d} m\right)\\
&=C\sum_{n=0}^{+\infty}\frac{(1-b)^p}{b^{p}}\int_{{B_{n+1}\backslash\overline{B_n}}}b^{p(n+1)}\mathrm{d}\Gamma(f)+C\sum_{n=0}^{+\infty}(1-b)^p\frac{b^{pn}}{\Psi(s_{n+1})}\int_{{B_{n+1}\backslash\overline{B_n}}}|f|^p\mathrm{d} m,
\end{align*}
where
$$\int_{{B_{n+1}\backslash\overline{B_n}}}b^{p(n+1)}\mathrm{d}\Gamma(f)\le\int_{{B_{n+1}\backslash\overline{B_n}}}|\widetilde{\phi}|^p\mathrm{d}\Gamma(f),$$
and
$$\frac{b^{pn}}{\Psi(s_{n+1})}=\frac{b^{pn}}{\Psi(r)}\frac{\Psi(r)}{\Psi(s_{n+1})}\le \frac{b^{pn}}{\Psi(r)}C_{\Psi}\left(\frac{r}{s_{n+1}}\right)^{\log_2C_\Psi}=C_\Psi \frac{b^{pn}}{\Psi(r)}\left(\frac{1}{(1-a)a^{n}}\right)^{\log_2C_\Psi}.$$
Let $b^p=a^{{\log_2C_\Psi}}$, then
\begin{align*}
&\int_{B(x,R+r)\backslash\overline{B(x,R)}}|\widetilde{f}|^p\mathrm{d}\Gamma(\phi)\\
&\le\frac{C(1-b)^p}{b^{p}}\int_{B(x,R+r)\backslash\overline{B(x,R)}}|\widetilde{\phi}|^p\mathrm{d}\Gamma(f)\\
&\hspace{20pt}+\frac{CC_\Psi(1-b)^p}{(1-b^{\frac{p}{\log_2C_\Psi}})^{\log_2C_\Psi}}\frac{1}{\Psi(r)}\int_{B(x,R+r)\backslash\overline{B(x,R)}}|f|^p\mathrm{d} m.
\end{align*}
For any $\delta>0$, since $\lim_{b\uparrow1}\frac{(1-b)^p}{b^{p}}=0$, there exists $b\in(0,1)$ depending on $C, \delta$, such that $\frac{C(1-b)^p}{b^{p}}<\delta$, let $C_\delta=\frac{CC_\Psi(1-b)^p}{(1-b^{\frac{p}{\log_2C_\Psi}})^{\log_2C_\Psi}}$, then we have (\ref{eq_CS_self}). Moreover, by taking $f\in\mathcal{F}\cap C_c(X)$ with $f\equiv1$ in $B(x,R+r)$, we have $\phi\in\mathcal{F}$.
\end{proof}

Let $U$ be an open subset of $X$. Let
$$\mathcal{F}(U)=\text{the }\mathcal{E}_1\text{-closure of }\mathcal{F}\cap C_c(U).$$
We say that $u\in\mathcal{F}$ is harmonic in $U$ if $\mathcal{E}(u;v)=0$ for any $v\in\mathcal{F}\cap C_c(U)$. We have the following characterization of harmonic functions.

\begin{lemma}\label{lem_harm_equiv}
Let $U$ be a bounded open set and $u\in\mathcal{F}$. The followings are equivalent.
\begin{enumerate}[label=(\arabic*)]
\item $u$ is harmonic in $U$, that is, $\mathcal{E}(u;v)=0$ for any $v\in\mathcal{F}\cap C_c(U)$.
\item $\mathcal{E}(u;v)=0$ for any $v\in\mathcal{F}$ with $\widetilde{v}=0$ q.e. on $X\backslash U$.
\item $\mathcal{E}(u)=\inf\{\mathcal{E}(v):v\in\mathcal{F},\widetilde{v}=\widetilde{u}\text{ q.e. on }X\backslash U\}$.
\end{enumerate}
\end{lemma}

\begin{proof}
``(2)$\Rightarrow$(1)": Trivial. ``(1)$\Rightarrow$(2)": It is obvious by Lemma \ref{lem_calFU}.

``(3)$\Rightarrow$(2)": For any $v\in\mathcal{F}$ with $\widetilde{v}=0$ q.e. on $X\backslash U$, for any $t>0$, we have $\widetilde{u}\pm t\widetilde{v}=\widetilde{u}$ q.e. on $X\backslash U$. By assumption, we have $\mathcal{E}(u)\le\mathcal{E}(u\pm tv)$, hence $\mathcal{E}(u;v)=\frac{1}{p}\lim_{t\downarrow0}\frac{1}{t}(\mathcal{E}(u+tv)-\mathcal{E}(u))\ge0$, $\mathcal{E}(u;v)=\frac{1}{p}\lim_{t\downarrow0}\frac{1}{-t}(\mathcal{E}(u-tv)-\mathcal{E}(u))\le0$, which gives $\mathcal{E}(u;v)=0$.

``(2)$\Rightarrow$(3)": We only need to show ``$\le$". For any $v\in\mathcal{F}$ with $\widetilde{v}=\widetilde{u}$ q.e. on $X\backslash U$, let $\varphi(t)=\mathcal{E}((1-t)u+tv)=\mathcal{E}(u+t(v-u))$, $t\in\mathbb{R}$, then $\varphi'(t)=\mathcal{E}(u+t(v-u);v-u)$. Since $\widetilde{v}-\widetilde{u}=0$ q.e. on $X\backslash U$, by assumption, we have $\varphi'(0)=\mathcal{E}(u;v-u)=0$. By (\ref{eq_quasi_strict}), we have $\varphi'(t)\ge\varphi'(0)=0$ for any $t\ge0$, hence $\mathcal{E}(v)=\varphi(1)\ge\varphi(0)=\mathcal{E}(u)$.
\end{proof}

The following reverse Poincar\'e inequality follows from the cutoff Sobolev inequality; this result parallels \cite[LEMMA 3.3]{KM20}.

\begin{proposition}[Reverse Poincar\'e inequality]\label{prop_PI_reverse}
Assume \ref{eq_VD}, \ref{eq_CS}. Then there exists $C>0$ such that for any ball $B(x,r)$, for any $h\in\mathcal{F}\cap L^\infty(X;m)$ which is harmonic in $B(x,2r)$, for any $a\in\mathbb{R}$, we have
$$\int_{B(x,r)}\mathrm{d}\Gamma(h)\le \frac{C}{\Psi(r)}\int_{B(x,2r)\backslash\overline{B(x,r)}}|h-a|^p\mathrm{d} m.$$
\end{proposition}

\begin{proof}
For notational convenience, we may assume that all functions in $\mathcal{F}$ are quasi-cont-inuous. By replacing $h$ with $h-a\psi$, where $\psi\in\mathcal{F}\cap C_c(X)$ satisfies $0\le\psi\le1$ on $X$, $\psi=1$ in $B(x,2r)$, we may assume that $a=0$. Let $\varepsilon, \delta>0$ be chosen later. Let $\phi\in\mathcal{F}$ be a cutoff function for $B(x,r)\subseteq B(x,2r)$ given by Proposition \ref{prop_CS_self}, then $h\phi^p\in\mathcal{F}$ satisfies $h\phi^p=0$ q.e. on $X\backslash B(x,2r)$. By Lemma \ref{lem_harm_equiv}, we have
\begin{align*}
0=\mathcal{E}(h;h\phi^p)=\Gamma(h;h\phi^p)(X)=\int_{X}{\phi}^p\mathrm{d}\Gamma(h)+p\int_X{h}\phi^{p-1}\mathrm{d}\Gamma(h;\phi),
\end{align*}
where the last equality follows from the chain rule and the Leibniz rule for $\Gamma(\cdot;\cdot)$, see \cite[Theorem 4.15 (b), Theorem 5.12, Definition 4.14]{KS24a}. By \cite[Proposition 4.8]{KS24a} and Young's inequality, we have
\begin{align*}
&\lvert\int_X{h}\phi^{p-1}\mathrm{d}\Gamma(h;\phi)\rvert\\
&\le\left(\int_X\phi^{(p-1)\frac{p}{p-1}}\mathrm{d}\Gamma(h)\right)^{(p-1)/p}\left(\int_X|h|^p\mathrm{d}\Gamma(\phi)\right)^{1/p}\\
&\le\varepsilon\int_X\phi^p\mathrm{d}\Gamma(h)+C_\varepsilon\int_X|h|^p\mathrm{d}\Gamma(\phi),
\end{align*}
here $C_\varepsilon$ is some positive constant depending only on $p, \varepsilon$. By Proposition \ref{prop_CS_self}, we have
\begin{align*}
&\int_X|h|^p\mathrm{d}\Gamma(\phi)=\int_{B(x,2r)\backslash\overline{B(x,r)}}|h|^p\mathrm{d}\Gamma(\phi)\\
&\le\delta\int_{B(x,2r)\backslash\overline{B(x,r)}}|\phi|^p\mathrm{d}\Gamma(h)+\frac{C_\delta}{\Psi(r)}\int_{B(x,2r)\backslash\overline{B(x,r)}}|h|^p\mathrm{d} m,
\end{align*}
hence
\begin{align*}
&\lvert\int_X{h}\phi^{p-1}\mathrm{d}\Gamma(h;\phi)\rvert\\
&\le \varepsilon\int_X\phi^p\mathrm{d}\Gamma(h)+C_\varepsilon \delta\int_{B(x,2r)\backslash\overline{B(x,r)}}\phi^p\mathrm{d}\Gamma(h)+\frac{C_\varepsilon C_\delta }{\Psi(r)}\int_{B(x,2r)\backslash\overline{B(x,r)}}|h|^p\mathrm{d} m\\
&\le(\varepsilon+\delta C_\varepsilon)\int_X\phi^p\mathrm{d}\Gamma(h)+\frac{C_\varepsilon C_\delta }{\Psi(r)}\int_{B(x,2r)\backslash\overline{B(x,r)}}|h|^p\mathrm{d} m,
\end{align*}
which gives
\begin{align*}
0&\ge\int_X\phi^p\mathrm{d}\Gamma(h)-p(\varepsilon+\delta C_\varepsilon)\int_X\phi^p\mathrm{d}\Gamma(h)-\frac{pC_\varepsilon C_\delta }{\Psi(r)}\int_{B(x,2r)\backslash\overline{B(x,r)}}|h|^p\mathrm{d} m\\
&=\left(1-p\varepsilon-p\delta C_\varepsilon\right)\int_X\phi^p\mathrm{d}\Gamma(h)-\frac{pC_\varepsilon C_\delta }{\Psi(r)}\int_{B(x,2r)\backslash\overline{B(x,r)}}|h|^p\mathrm{d} m,
\end{align*}
that is,
$$\left(1-p\varepsilon-p\delta C_\varepsilon\right)\int_X\phi^p\mathrm{d}\Gamma(h)\le\frac{pC_\varepsilon C_\delta }{\Psi(r)}\int_{B(x,2r)\backslash\overline{B(x,r)}}|h|^p\mathrm{d} m.$$
Firstly, let $\varepsilon=\frac{1}{4p}$, secondly, let $\delta=\frac{1}{4pC_\varepsilon}$, then we have
$$\int_{B(x,r)}\mathrm{d}\Gamma(h)\le\int_X\phi^p\mathrm{d}\Gamma(h)\le\frac{2pC_\varepsilon C_\delta }{\Psi(r)}\int_{B(x,2r)\backslash\overline{B(x,r)}}|h|^p\mathrm{d} m.$$
\end{proof}

\section{Proof of Theorem \ref{thm_singular}}\label{sec_singular}

Since mutual singularity of measures is a local property, we may assume that $(X,d)$ is unbounded throughout this section.

We follow the strategy employed in \cite{KM20}. Firstly, we prove that any non-negative function in $\mathcal{F}\cap C_c(X)$ can be approximated by ``piecewise harmonic functions". We need do some preparations as follows.

\begin{lemma}\label{lem_calFU_Banach}
Assume \ref{eq_VD}, \ref{eq_PI}. Then the bottom spectrum positivity condition \ref{eq_BSP} holds: for any bounded open subset $U\subseteq X$,
\begin{equation*}\label{eq_BSP}\tag*{\text{BSP}}
\lambda_1(U)=\inf \left\{\frac{\mathcal{E}(u)}{\lVert u\rVert_{L^p(X;m)}^p}:u\in \mathcal{F}(U)\backslash \{0\}\right\}>0.
\end{equation*}
Equivalently, there exists $C>0$, which may depend on $U$, such that for any $u\in\mathcal{F}(U)$,
$$\int_U|u|^p\mathrm{d} m\le C\mathcal{E}(u).$$
\end{lemma}

\begin{proof}
Since $U$ is bounded, there exists a ball $B(x_0,r)\supseteq U$. On the one hand, by \ref{eq_PI},
\begin{align*}
&\int_{B(x_0,2r)}\int_{B(x_0,2r)}|u(x)-u(y)|^pm(\mathrm{d} x)m(\mathrm{d} y)\\
&\le2^{p-1}\int_{B(x_0,2r)}\int_{B(x_0,2r)}\left(|u(x)-u_{B(x_0,2r)}|^p+|u(y)-u_{B(x_0,2r)}|^p\right)m(\mathrm{d} x)m(\mathrm{d} y)\\
&=2^pm(B(x_0,2r))\int_{B(x_0,2r)}|u-u_{B(x_0,2r)}|^p\mathrm{d} m\\
&\le{2^pC_{PI}m(B(x_0,2r))}{\Psi(2r)}\int_{B(x_0,2A_{PI}r)}\mathrm{d}\Gamma(u)={2^pC_{PI}m(B(x_0,2r))}{\Psi(2r)}\mathcal{E}(u).
\end{align*}
On the other hand, by \ref{eq_CC}, there exists a ball $B(y_0,\frac{r}{4})\subseteq B(x_0,2r)\backslash B(x_0,r)$, hence
\begin{align*}
&\int_{B(x_0,2r)}\int_{B(x_0,2r)}|u(x)-u(y)|^pm(\mathrm{d} x)m(\mathrm{d} y)\\
&\ge\int_{U}\int_{B(y_0,\frac{r}{4})}|u(x)-u(y)|^pm(\mathrm{d} x)m(\mathrm{d} y)=m(B(y_0,\frac{r}{4}))\int_U|u|^p\mathrm{d} m.
\end{align*}
Therefore,
$$m(B(y_0,\frac{r}{4}))\int_U|u|^p\mathrm{d} m\le {2^pC_{PI}m(B(x_0,2r))}{\Psi(2r)}\mathcal{E}(u),$$
where by \ref{eq_VD},
$$\frac{m(B(x_0,2r))}{m(B(y_0,\frac{r}{4}))}\le C_{VD}\left(\frac{d(x_0,y_0)+2r}{\frac{r}{4}}\right)^{\log_2C_{VD}}\le C_{VD}^5,$$
which gives
$$\int_{U}|u|^p\mathrm{d} m\le {2^pC_{VD}^5C_{PI}}{\Psi(2r)} \mathcal{E}(u)\le {2^pC_{\Psi}C_{VD}^5C_{PI}}{\Psi(r)} \mathcal{E}(u).$$
\end{proof}

\begin{proposition}\label{prop_harm_exist}
Assume \ref{eq_VD}, \ref{eq_BSP}. Let $U$ be a bounded open set and $u\in\mathcal{F}$. Then there exists a unique function $h\in\mathcal{F}$ such that $h$ is harmonic in $U$ and $\widetilde{h}=\widetilde{u}$ q.e. on $X\backslash U$. We denote this function by $H^Uu$. Moreover, if $0\le u\le M$ $m$-a.e. in $X$, where $M>0$ is some constant, then $0\le H^Uu\le M$ $m$-a.e. in $X$.
\end{proposition}

\begin{proof}
Let
$$\mathcal{L}_{u,X\backslash U}=\left\{v\in\mathcal{F}:\widetilde{v}=\widetilde{u}\text{ q.e. on }X\backslash U\right\},$$
then by Lemma \ref{lem_calFU}, we have $\mathcal{L}_{u,X\backslash U}=\left\{u+v:v\in\mathcal{F}(U)\right\}$. Consider the variational problem $\inf\{\mathcal{E}(v):v\in\mathcal{L}_{u,X\backslash U}\}$. For any minimizing sequence $\{v_n\}\subseteq\mathcal{L}_{u,X\backslash U}$, since $\frac{v_n+v_m}{2}\in\mathcal{L}_{u,X\backslash U}$ for any $n,m$, by \ref{eq_Cla}, we have $\{v_n\}$ is an $\mathcal{E}$-Cauchy sequence. Since $v_n-v_m\in\mathcal{F}(U)$ for any $n,m$, by Lemma \ref{lem_calFU_Banach}, we have $\{v_n\}$ is an $\mathcal{E}_1$-Cauchy sequence, then there exists $h\in\mathcal{F}$ such that $\{v_n\}$ is $\mathcal{E}_1$-convergent to $h$. By \cite[Corollary 8.7]{Yan25a}, there exists a subsequence, still denoted by $\{v_n\}$, such that $\{\widetilde{v}_n\}$ converges to $\widetilde{h}$ q.e. on $X$, which gives $\widetilde{h}=\widetilde{u}$ q.e. on $X\backslash U$, that is, $h\in\mathcal{L}_{u,X\backslash U}$. Hence $h\in\mathcal{F}$ satisfies $\mathcal{E}(h)=\inf\{\mathcal{E}(v):v\in\mathcal{F},\widetilde{v}=\widetilde{u}\text{ q.e. on }X\backslash U\}=\inf\{\mathcal{E}(v):v\in\mathcal{F},\widetilde{v}=\widetilde{h}\text{ q.e. on }X\backslash U\}$, by Lemma \ref{lem_harm_equiv}, we have $h$ is harmonic in $U$.

Let $g\in\mathcal{F}$ satisfy that $g$ is harmonic in $U$ and $\widetilde{g}=\widetilde{u}$ q.e. on $X\backslash U$. Let $\varphi(t)=\mathcal{E}((1-t)h+tg)=\mathcal{E}(h+t(g-h))$, $t\in\mathbb{R}$, then $\varphi'(t)=\mathcal{E}(h+t(g-h);g-h)$. Since $g,h$ are harmonic in $U$ and $\widetilde{g}=\widetilde{h}=\widetilde{u}$ q.e. on $X\backslash U$, by Lemma \ref{lem_harm_equiv}, we have $\varphi'(0)=\mathcal{E}(h;g-h)=0$, $\varphi'(1)=\mathcal{E}(g;g-h)=0$, by (\ref{eq_quasi_strict}), we have $g=h$.

Moreover, if $0\le u\le M$ $m$-a.e. in $X$, then $(\widetilde{h}\vee0)\wedge M=(\widetilde{u}\vee0)\wedge M=\widetilde{u}$ q.e. on $X\backslash U$. Since $\mathcal{E}((h\vee 0)\wedge M)\le\mathcal{E}(h)=\inf\{\mathcal{E}(v):v\in\mathcal{F},\widetilde{v}=\widetilde{u}\text{ q.e. on }X\backslash U\}=\inf\{\mathcal{E}(v):v\in\mathcal{F},\widetilde{v}=(\widetilde{h}\vee0)\wedge M\text{ q.e. on }X\backslash U\}$, by Lemma \ref{lem_harm_equiv}, we have $(h\vee0)\wedge M\in\mathcal{F}$ is harmonic in $U$. By the above uniqueness result, $(h\vee0)\wedge M=h$, that is, $0\le h\le M$ $m$-a.e. in $X$.
\end{proof}

Any non-negative function in $\mathcal{F}\cap C_c(X)$ can be approximated by ``piecewise harmonic functions" as follows; this result parallels \cite[PROPOSITION 3.9]{KM20}.

\begin{proposition}\label{prop_harm_approx}
Assume \ref{eq_VD}, \ref{eq_BSP}. For any $f\in\mathcal{F}\cap C_c(X)$ with $0\le f\le 1$ on $X$, for any $n\ge0$, let $F_n=\{x\in X:f(x)=\frac{k}{2^n},k=0,1,\ldots,2^n\}$, for any $k=0,1,\ldots,2^n-1$, let $U_{n,k}=\{x\in X:\frac{k}{2^n}<f(x)<\frac{k+1}{2^n}\}$ and $f_{n,k}=H^{U_{n,k}}(((f-\frac{k}{2^n})\vee0)\wedge \frac{1}{2^n})$, let $f_n=\sum_{k=0}^{2^n-1}f_{n,k}$. Then for any $n\ge0$, $\widetilde{f}_n=f$ q.e on $F_n$, $f_n$ is harmonic in $X\backslash F_n$, $\Gamma(f_n)(F_n)=0$ and $|f_n-f|\le \frac{1}{2^n}1_{\{f>0\}}$. Moreover, $\lim_{n\to+\infty}\mathcal{E}_1(f_n-f)=0$.
\end{proposition}

\begin{proof}
For notational convenience, we may assume that all functions in $\mathcal{F}$ are quasi-cont-inuous. Since $U_{n,k}$ is a bounded open set, by Proposition \ref{prop_harm_exist}, $f_{n,k}\in\mathcal{F}$ is harmonic in $U_{n,k}$, $0\le f_{n,k}\le \frac{1}{2^n}$ in $X$, $f_{n,k}=\frac{1}{2^n}$ q.e. on $\{f\ge \frac{k+1}{2^n}\}$, $f_{n,k}=0$ q.e. on $\{f\le \frac{k}{2^n}\}$, which gives $f_n\in\mathcal{F}$ and $0\le f_n\le 1$ in $X$. In particular, $f_{n,k}$ is harmonic in $\{f\not\in\{\frac{k}{2^n},\frac{k+1}{2^n}\}\}$, since for any $v\in\mathcal{F}\cap C_c(\{f\not\in\{\frac{k}{2^n},\frac{k+1}{2^n}\}\})$, $v=v1_{U_{n,k}}+v1_{\{f<\frac{k}{2^n}\}}+v1_{\{f>\frac{k+1}{2^n}\}}$, where $v1_{U_{n,k}}\in\mathcal{F}\cap C_c(U_{n,k})$, $v1_{\{f<\frac{k}{2^n}\}}\in\mathcal{F}\cap C_c(\{f<\frac{k}{2^n}\})$, $v1_{\{f>\frac{k+1}{2^n}\}}\in\mathcal{F}\cap C_c(\{f>\frac{k+1}{2^n}\})$, which gives $\mathcal{E}(f_{n,k};v)=\mathcal{E}(f_{n,k};v1_{U_{n,k}})+\mathcal{E}(f_{n,k};v1_{\{f<\frac{k}{2^n}\}})+\mathcal{E}(f_{n,k};v1_{\{f>\frac{k+1}{2^n}\}})=0$. Hence $f_n$ is harmonic in $X\backslash F_n$. For any $k=0,1,\ldots,2^n$, $f_{n,l}=\frac{1}{2^n}$ q.e. on $\{f=\frac{k}{2^n}\}$ if $l<k$, and $f_{n,l}=0$ q.e. on $\{f=\frac{k}{2^n}\}$ if $l\ge k$, hence $f_n=\sum_{l=0}^{k-1}f_{n,l}=\frac{k}{2^n}=f$ q.e. on $\{f=\frac{k}{2^n}\}$, that is, $f_n=f$ q.e. on $F_n$. By \cite[Proposition 8.12]{Yan25a}, $\Gamma(f_n)(F_n\cap\{f_n\ne f\})=0$.

By \cite[Theorem 4.17, Theorem 4.18]{KS24a} or an argument similar to that of \cite[THEOREM 4.3.8]{CF12}, we have $\Gamma(f_n)(F_n)\le\Gamma(f_n)(f_n^{-1}(\{\frac{k}{2^n}:k=0,1,\ldots,2^n\}))=0$ and $\Gamma(f)(F_n)=0$.

For any $k=0,1,\ldots,2^n-1$, we have $f_{n,l}=\frac{1}{2^n}$ q.e. on $U_{n,k}$ if $l<k$, and $f_{n,l}=0$ q.e. on $U_{n,k}$ if $l>k$, hence $f_n=\sum_{l=0}^kf_{n,l}=\frac{k}{2^n}+f_{n,k}\in[\frac{k}{2^n},\frac{k+1}{2^n}]$ q.e. on $U_{n,k}$. Since $\frac{k}{2^n}<f< \frac{k+1}{2^n}$ on $U_{n,k}$, we have $|f_n-f|\le \frac{1}{2^n}$ q.e. on $U_{n,k}$. Since $f_n=f$ q.e. on $F_n$, we have $|f_n-f|\le \frac{1}{2^n}1_{\{f>0\}}$. By the dominated convergence theorem, we have $\lim_{n\to+\infty} \lVert {f_n-f}\rVert_{L^p(X;m)}=0$.

For any $n\ge0$, by the strongly local property and Lemma \ref{lem_harm_equiv}, we have
$$\mathcal{E}(f_n)=\sum_{k=0}^{2^n-1}\Gamma(f_n)(U_{n,k})=\sum_{k=0}^{2^n-1}\Gamma(f_{n,k})(U_{n,k})\le\sum_{k=0}^{2^n-1}\Gamma(f)(U_{n,k})\le\mathcal{E}(f),$$
and
$$\mathcal{E}(f_n)=\sum_{k=0}^{2^n-1}\Gamma(f_{n,k})(U_{n,k})\le\sum_{k=0}^{2^n-1}\Gamma(f_{n+1})(U_{n,k})\le\mathcal{E}(f_{n+1}),$$
hence $\{\mathcal{E}(f_n)\}$ is an increasing sequence bounded above by $\mathcal{E}(f)$, which implies that $A=\lim_{n\to+\infty}\mathcal{E}(f_n)\le\mathcal{E}(f)$.

We claim that $\{f_n\}$ is an $\mathcal{E}$-Cauchy sequence. Indeed, for any $n\le m$, by the strongly local property and Lemma \ref{lem_harm_equiv}, we have
$$\mathcal{E}(f_n)=\sum_{k=0}^{2^n-1}\Gamma(f_{n,k})(U_{n,k})\le\sum_{k=0}^{2^n-1}\Gamma(\frac{f_n+f_m}{2})(U_{n,k})\le\mathcal{E}(\frac{f_n+f_m}{2}).$$
If $p\in(1,2]$, then by \ref{eq_Cla}, we have
$$2 \left(\mathcal{E}(\frac{f_n+f_m}{2})^{\frac{1}{p-1}}+\mathcal{E}(\frac{f_n-f_m}{2})^{\frac{1}{p-1}}\right)^{p-1}\le\mathcal{E}(f_n)+\mathcal{E}(f_m)\to 2A,$$
where
$$\text{LHS}\ge 2 \left(\mathcal{E}(f_n)^{\frac{1}{p-1}}+\mathcal{E}(\frac{f_n-f_m}{2})^{\frac{1}{p-1}}\right)^{p-1}\to2 \left(A^{\frac{1}{p-1}}+\varlimsup_{n\le m\to+\infty}\mathcal{E}(\frac{f_n-f_m}{2})^{\frac{1}{p-1}}\right)^{p-1},$$
hence $\lim_{n\le m\to+\infty}\mathcal{E}(\frac{f_n-f_m}{2})=0$. If $p\in[2,+\infty)$, then by \ref{eq_Cla}, we have
$$\mathcal{E}(\frac{f_n+f_m}{2})+\mathcal{E}(\frac{f_n-f_m}{2})\le2 \left(\mathcal{E}(\frac{f_n}{2})^{\frac{1}{p-1}}+\mathcal{E}(\frac{f_m}{2})^{\frac{1}{p-1}}\right)^{p-1}\to A,$$
where
$$\text{LHS}\ge\mathcal{E}(f_n)+\mathcal{E}(\frac{f_n-f_m}{2})\to A+\varlimsup_{n\le m\to+\infty}\mathcal{E}(\frac{f_n-f_m}{2}),$$
hence $\lim_{n\le m\to+\infty}\mathcal{E}(\frac{f_n-f_m}{2})=0$. Therefore, $\{f_n\}$ is an $\mathcal{E}$-Cauchy sequence. Since $\{f_n\}$ is $L^p$-convergent to $f$, we have $\{f_n\}$ is $\mathcal{E}_1$-convergent to $f$.
\end{proof}

Secondly, we prove the singularity result for bounded harmonic functions as follows; this result parallels \cite[PROPOSITION 3.5]{KM20}.

\begin{proposition}\label{prop_har_singular}
Assume \ref{eq_VD}, \ref{eq_PI}, \ref{eq_CS}, (\ref{eq_Psi_singular}), and $d$ is geodesic. Let $U$ be an open set and $h\in\mathcal{F}\cap L^\infty(X;m)$ harmonic in $U$. Then $\Gamma(h)|_U\perp m|_U$.
\end{proposition}

\begin{proof}
Suppose that the result is not true. Let $A=A_{PI}$ be the constant in \ref{eq_PI}. By \cite[LEMMA 3.1]{KM20}, for any $\lambda>4A$, there exist $x\in U$, $r_0=r_0(x,\lambda)>0$ with $B(x,r_0)\subseteq U$ such that for any $r\in(0,r_0)$, for any $\delta\in[\frac{1}{\lambda},1]$, for any $y\in B(x,r)$, we have
\begin{equation}\label{eq_har_singular1}
0<\frac{1}{2}\frac{\mathrm{d}\Gamma(h)}{\mathrm{d} m}(x)\le \frac{\Gamma(h)(B(y,\delta r))}{m(B(y,\delta r))}\le 2\frac{\mathrm{d}\Gamma(h)}{\mathrm{d} m}(x).
\end{equation}

Let $V$ be an $\frac{r}{\lambda}$-net in $B(x,r)$. We claim that for any $y_1$, $y_2\in V$, we have
\begin{equation}\label{eq_har_singular2}
|h_{B(y_1,\frac{r}{\lambda})}-h_{B(y_2,\frac{r}{\lambda})}|\lesssim\lambda\Psi(\frac{r}{\lambda})^{1/p}.
\end{equation}
Firstly, for any $y_1$, $y_2\in V$ with $d(y_1,y_2)<\frac{3r}{\lambda}$, we have
\begin{align*}
&|h_{B(y_1,\frac{r}{\lambda})}-h_{B(y_2,\frac{r}{\lambda})}|^p\le\dashint_{B(y_1,\frac{r}{\lambda})}\dashint_{B(y_2,\frac{r}{\lambda})}|h(x_1)-h(x_2)|^pm(\mathrm{d} x_1)m(\mathrm{d} x_2)\\
&\le \frac{1}{m(B(y_1,\frac{r}{\lambda}))m(B(y_2,\frac{r}{\lambda}))}\int_{B(y_1,\frac{4r}{\lambda})}\int_{B(y_1,\frac{4r}{\lambda})}|h(x_1)-h(x_2)|^pm(\mathrm{d} x_1)m(\mathrm{d} x_2)\\
&\le\frac{2^pm(B(y_1,\frac{4r}{\lambda}))}{m(B(y_1,\frac{r}{\lambda}))m(B(y_2,\frac{r}{\lambda}))}\int_{B(y_1,\frac{4r}{\lambda})}|h-h_{B(y_1,\frac{4r}{\lambda})}|^p\mathrm{d} m\\
&\overset{(\dagger)}{\scalebox{2}[1]{$\lesssim$}}\frac{\Psi(\frac{4r}{\lambda})}{m(B(y_1,\frac{4Ar}{\lambda}))}\int_{B(y_1,\frac{4Ar}{\lambda})}\mathrm{d}\Gamma(h)\overset{(\diamond)}{\scalebox{2}[1]{$\lesssim$}}\Psi(\frac{r}{\lambda}),
\end{align*}
where $(\dagger)$ follows from \ref{eq_PI} and \ref{eq_VD}, and $(\diamond)$ follows from (\ref{eq_har_singular1}) and the doubling property of $\Psi$.

Secondly, for any $y_1$, $y_2\in V$, there exists a sequence $\{z_k:k=0,1,\ldots, K\}\subseteq V$ with $K\le 3\lambda$ such that $z_0=y_1$, $z_K=y_2$ and $d(z_k,z_{k+1})<\frac{3r}{\lambda}$ for any $k=0,1,\ldots,K-1$. Indeed, let $\gamma$ be a geodesic in $B(x,r)$ connecting $x$, $y_1$, let $K_1\ge0$ be the integer satisfying $K_1\le \frac{d(x,y_1)}{\frac{r}{\lambda}}<K_1+1$, then $K_1<\lambda$, there exist ${x}_0$, \ldots, ${x}_{K_1}$, ${x}_{K_1+1}\in\gamma$ with $x_0=x$, $x_{K_1+1}=y_1$ such that $d(x_k,x_{k+1})\le \frac{r}{\lambda}$ for any $k=0,1,\ldots,K_1$. For any $k=0,1,\ldots,K_1+1$, there exists $\widetilde{z}_k\in V$ such that $d(\widetilde{z}_k,x_k)<\frac{r}{\lambda}$, in particular, $\widetilde{z}_{K_1+1}=y_1$. Hence $d(\widetilde{z}_0,x)<\frac{r}{\lambda}$, $\widetilde{z}_{K_1+1}=y_1$ and $d(\widetilde{z}_k,\widetilde{z}_{k+1})<\frac{3r}{\lambda}$ for any $k=0,1,\ldots,K_1$. Similarly, there exist $\widehat{z}_0$, \ldots, $\widehat{z}_{K_2}$, $\widehat{z}_{K_2+1}$ with $0\le K_2<\lambda$ such that $d(\widehat{z}_0,x)<\frac{r}{\lambda}$, $\widehat{z}_{K_2+1}=y_2$ and $d(\widehat{z}_k,\widehat{z}_{k+1})<\frac{3r}{\lambda}$ for any $k=0,1,\ldots,K_2$. By concatenating $\{\widetilde{z}_k\}$ and $\{\widehat{z}_k\}$, we obtain $\{z_k\}$. Therefore
$$|h_{B(y_1,\frac{r}{\lambda})}-h_{B(y_2,\frac{r}{\lambda})}|\le\sum_{k=0}^{K-1}|h_{B(z_k,\frac{r}{\lambda})}-h_{B(z_{k+1},\frac{r}{\lambda})}|\lesssim \lambda\Psi(\frac{r}{\lambda})^{1/p},$$
that is, the claim holds true.

For any $c\in\mathbb{R}$, we have
$$\int_{B(x,r)}|h-h_{B(x,r)}|^p\mathrm{d} m\le2^{p-1}\int_{B(x,r)}\left(|h-c|^p+|h_{B(x,r)}-c|^p\right)\mathrm{d} m,$$
where
$$|h_{B(x,r)}-c|^p\le\left(\dashint_{B(x,r)}|h-c|\mathrm{d} m\right)^p\le\dashint_{B(x,r)}|h-c|^p\mathrm{d} m,$$
hence
$$\int_{B(x,r)}|h-h_{B(x,r)}|^p\mathrm{d} m\le2^p\int_{B(x,r)}|h-c|^p\mathrm{d} m\text{ for any }c\in\mathbb{R}.$$
Take $y_1\in V$, let $c=h_{B(y_1,\frac{r}{\lambda})}$, then
\begin{align*}
&\int_{B(x,r)}|h-h_{B(x,r)}|^p\mathrm{d} m\le2^p\int_{B(x,r)}|h-h_{B(y_1,\frac{r}{\lambda})}|^p\mathrm{d} m\\
&\le2^p\sum_{y_2\in V}\int_{B(y_2,\frac{r}{\lambda})}|h-h_{B(y_1,\frac{r}{\lambda})}|^p\mathrm{d} m\\
&\le2^{2p-1}\sum_{y_2\in V}\int_{B(y_2,\frac{r}{\lambda})}\left(|h-h_{B(y_2,\frac{r}{\lambda})}|^p+|h_{B(y_1,\frac{r}{\lambda})}-h_{B(y_2,\frac{r}{\lambda})}|^p\right)\mathrm{d} m\\
&\overset{(\star)}{\scalebox{2}[1]{$\lesssim$}}\sum_{y_2\in V}\Psi(\frac{r}{\lambda})\int_{B(y_2,\frac{Ar}{\lambda})}\mathrm{d}\Gamma(h)+\sum_{y_2\in V}\lambda^p\Psi(\frac{r}{\lambda})m(B(y_2,\frac{r}{\lambda}))\\
&\overset{(\diamond)}{\scalebox{2}[1]{$\lesssim$}}\sum_{y_2\in V}\Psi(\frac{r}{\lambda})m(B(y_2,\frac{Ar}{\lambda}))+\lambda^p\Psi(\frac{r}{\lambda})m(B(x,r))\\
&\overset{(\dagger)}{\scalebox{2}[1]{$\lesssim$}}\Psi(\frac{r}{\lambda})m(B(x,r))+\lambda^p\Psi(\frac{r}{\lambda})m(B(x,r))\asymp\lambda^p\Psi(\frac{r}{\lambda})m(B(x,r)),
\end{align*}
where ($\star$) follows from \ref{eq_PI} and (\ref{eq_har_singular2}), ($\diamond$) follows from (\ref{eq_har_singular1}) and \ref{eq_VD}, and ($\dagger$) follows from \ref{eq_VD}. However, by Proposition \ref{prop_PI_reverse}, we have
$$\int_{B(x,r)}|h-h_{B(x,r)}|^p\mathrm{d} m\gtrsim\Psi(r)\Gamma(h)(B(x,\frac{r}{2}))\gtrsim\Psi(r)m(B(x,r)),$$
where the last inequality follows from (\ref{eq_har_singular1}) and \ref{eq_VD}. Hence $\Psi(r)\lesssim\lambda^p\Psi(\frac{r}{\lambda})$, that is, $\frac{\lambda^p\Psi(\frac{r}{\lambda})}{\Psi(r)}\gtrsim1$ for any $r\in(0,r_0(x,\lambda))$, for any $\lambda>4A$, which gives $\varliminf_{\lambda\to+\infty}\varliminf_{r\downarrow0}\frac{\lambda^p\Psi(\frac{r}{\lambda})}{\Psi(r)}>0$, contradicting to (\ref{eq_Psi_singular}).
\end{proof}

Thirdly, the singularity is preserved under linear combinations and $\mathcal{E}$-convergence; the proof is along the same lines as that of \cite[LEMMA 3.6 (b)]{KM20} and \cite[LEMMA 3.7 (b)]{KM20}.

\begin{lemma}\label{lem_sing_prep}
\hspace{0em}
\begin{enumerate}[label=(\arabic*),ref=(\arabic*)]
\item\label{lem_sing_prep1} If $f, g\in\mathcal{F}$ satisfy that $\Gamma(f)\perp m$ and $\Gamma(g)\perp m$, then for any $a,b\in\mathbb{R}$, we have $\Gamma(af+bg)\perp m$.
\item\label{lem_sing_prep2} If $\{f_n\}\subseteq\mathcal{F}$ and $f\in\mathcal{F}$ satisfy that $\Gamma(f_n)\perp m$ for any $n$, and $\lim_{n\to+\infty}\mathcal{E}(f_n-f)=0$, then $\Gamma(f)\perp m$.
\end{enumerate} 
\end{lemma}

Finally, we prove the singularity result for all functions in $\mathcal{F}$.

\begin{proof}[Proof of Theorem \ref{thm_singular}]
By \ref{eq_CC}, we have $d$ is Lipschitz equivalent to a geodesic metric, see \cite[PROPOSITION A.1]{KM20}. Since \ref{eq_VD}, \ref{eq_PI}, \ref{eq_CS} are stable under Lipschitz equivalence, we may assume that $d$ is geodesic. By Lemma \ref{lem_sing_prep}, we only need to prove for $f\in\mathcal{F}\cap C_c(X)$ with $0\le f\le1$ on $X$. Let $f_n, F_n$ be given by Proposition \ref{prop_harm_approx}. Since $f_n$ is harmonic in $X\backslash F_n$, by Proposition \ref{prop_har_singular}, we have $\Gamma(f_n)|_{X\backslash F_n}\perp m|_{X\backslash F_n}$. Since $\Gamma(f_n)(F_n)=0$, we have $\Gamma(f_n)\perp m$. Since $\lim_{n\to+\infty}\mathcal{E}_1(f_n-f)=0$, by Lemma \ref{lem_sing_prep} \ref{lem_sing_prep2}, we have $\Gamma(f)\perp m$.
\end{proof}

\section{Proof of Theorem \ref{thm_equiv}}\label{sec_equiv}

The key ingredient in the proof is the following result.

\begin{proposition}\label{prop_CS}
Assume \ref{eq_SVR}, \ref{eq_PI}, \ref{eq_cap}. Then \ref{eq_CS} holds.
\end{proposition}

We use an idea introduced by the author in \cite{Yan25b}.

\begin{lemma}[Morrey-Sobolev inequality]\label{lem_MS}
Assume \ref{eq_SVR}, \ref{eq_PI}. Then there exists $C_{MS}>0$ such that for any $u\in\mathcal{F}$, for any ball $B(x_0,R)$, for $m$-a.e. $x,y\in B(x_0,R)$ with $x\ne y$, we have
$$|u(x)-u(y)|^p\le C_{MS} \frac{\Psi(d(x,y))}{\Phi(d(x,y))}\int_{B(x_0,8A_{PI}R)}\mathrm{d}\Gamma(u),$$
where $A_{PI}$ is the constant in \ref{eq_PI}. Hence any function in $\mathcal{F}$ has a continuous version, or equivalently, $\mathcal{F}\subseteq C(X)$.
\end{lemma}

\begin{proof}
The proof is standard using the telescopic technique, see also \cite[Page 1654]{BCK05}. Let $x,y$ be two different Lebesgue points of $u\in\mathcal{F}\subseteq L^p(X;m)$, denote $r=d(x,y)$, then
$$|u(x)-u(y)|\le|u(x)-u_{B(x,r)}|+|u_{B(x,r)}-u_{B(y,r)}|+|u(y)-u_{B(y,r)}|.$$
Let $A=2A_{PI}\ge2$. For any $n\in\mathbb{Z}$, denote $B_n=B(x,A^{-n}r)$. For any $n\ge0$, by \ref{eq_VPhi} and \ref{eq_PI}, we have
\begin{align*}
&|u_{B_{n+1}}-u_{B_n}|\le\dashint_{B_{n+1}}|u-u_{B_n}|\mathrm{d} m\lesssim\dashint_{B_n}|u-u_{B_n}|\mathrm{d} m\le\left(\dashint_{B_n}|u-u_{B_n}|^p\mathrm{d} m\right)^{1/p}\\
&\lesssim\left(\frac{1}{\Phi(A^{-n}r)}\Psi(A^{-n}r)\int_{A_{PI}B_{n}}\mathrm{d}\Gamma(u)\right)^{1/p}\le\left(\frac{\Psi(A^{-n}r)}{\Phi(A^{-n}r)}\int_{B(x_0,8A_{PI}R)}\mathrm{d}\Gamma(u)\right)^{1/p},
\end{align*}
hence
\begin{align*}
&|u(x)-u_{B(x,r)}|=|u(x)-u_{B_0}|\\
&=\lim_{n\to+\infty}|u_{B_{n}}-u_{B_0}|\le\sum_{n=0}^{+\infty}|u_{B_{n+1}}-u_{B_n}|\\
&\lesssim\sum_{n=0}^{+\infty}\left(\frac{\Psi(A^{-n}r)}{\Phi(A^{-n}r)}\right)^{1/p}\left(\int_{B(x_0,8A_{PI}R)}\mathrm{d}\Gamma(u)\right)^{1/p}.
\end{align*}
By \ref{eq_SVR}, we have
\begin{align*}
&\sum_{n=0}^{+\infty}\left(\frac{\Psi(A^{-n}r)}{\Phi(A^{-n}r)}\right)^{1/p}\le\sum_{n=0}^{+\infty}\left(\Theta(A^{-n})\frac{\Psi(r)}{\Phi(r)}\right)^{1/p}\le\left(\frac{\Psi(r)}{\Phi(r)}\right)^{1/p}\sum_{n=0}^{+\infty}\Theta\left(2^{-n}\right)^{1/p},
\end{align*}
hence
$$|u(x)-u_{B(x,r)}|\lesssim\left(\frac{\Psi(r)}{\Phi(r)}\int_{B(x_0,8A_{PI}R)}\mathrm{d}\Gamma(u)\right)^{1/p}.$$
Similarly, we have
$$|u(y)-u_{B(y,r)}|\lesssim\left(\frac{\Psi(r)}{\Phi(r)}\int_{B(x_0,8A_{PI}R)}\mathrm{d}\Gamma(u)\right)^{1/p}.$$
By \ref{eq_VPhi} and \ref{eq_PI}, we have
\begin{align*}
&|u_{B(x,r)}-u_{B(y,r)}|\le\dashint_{B(x,r)}\dashint_{B(y,r)}|u(z_1)-u(z_2)|m(\mathrm{d} z_1)m(\mathrm{d} z_2)\\
&\le\left(\dashint_{B(x,r)}\dashint_{B(y,r)}|u(z_1)-u(z_2)|^pm(\mathrm{d} z_1)m(\mathrm{d} z_2)\right)^{1/p}\\
&\lesssim\left(\frac{1}{\Phi(r)^2}\int_{B(x,2r)}\int_{B(x,2r)}|u(z_1)-u(z_2)|^pm(\mathrm{d} z_1)m(\mathrm{d} z_2)\right)^{1/p}\\
&\le\left(\frac{2^pm(B(x,2r))}{\Phi(r)^2}\int_{B(x,2r)}|u-u_{B(x,2r)}|^p\mathrm{d} m\right)^{1/p}\\
&\lesssim\left(\frac{\Psi(r)}{\Phi(r)}\int_{B(x,2A_{PI}r)}\mathrm{d}\Gamma(u)\right)^{1/p}\le\left(\frac{\Psi(r)}{\Phi(r)}\int_{B(x_0,8A_{PI}R)}\mathrm{d}\Gamma(u)\right)^{1/p}.
\end{align*}
In summary, we have
$$|u(x)-u(y)|\lesssim\left(\frac{\Psi(d(x,y))}{\Phi(d(x,y))}\int_{B(x_0,8A_{PI}R)}\mathrm{d}\Gamma(u)\right)^{1/p}.$$
By \ref{eq_SVR}, we have
$$\frac{\Psi(d(x,y))}{\Phi(d(x,y))}\le\Theta \left(\frac{d(x,y)}{2R}\right)\frac{\Psi(2R)}{\Phi(2R)}.$$
Since $\lim_{t\downarrow0}\Theta(t)=\Theta(0)=0$, we have $\lim_{d(x,y)\downarrow0}|u(x)-u(y)|=0$, which implies that $u$ has a continuous version.
\end{proof}

\begin{proof}[Proof of Proposition \ref{prop_CS}]
For any ball $B(x_0,R)$, by \ref{eq_cap}, there exists a cutoff function $\phi\in\mathcal{F}$ for $B(x_0,R)\subseteq B(x_0,A_{cap}R)$ such that
$$\mathcal{E}(\phi)=\int_{B(x_0,A_{cap}R)}\mathrm{d}\Gamma(\phi)\le2C_{cap}\frac{V(x_0,R)}{\Psi(R)}.$$
By Lemma \ref{lem_MS}, for any $f\in\mathcal{F}\subseteq C(X)$, for any $x,y\in B(x_0,A_{cap}R)$, by \ref{eq_SVR}, we have
\begin{align*}
&|f(x)-f(y)|^p\le C_{MS}\frac{\Psi(d(x,y))}{\Phi(d(x,y))}\int_{B(x_0,8A_{PI}A_{cap}R)}\mathrm{d}\Gamma(f)\\
&\le C_{MS}\Theta(1)\frac{\Psi(2A_{cap}R)}{\Phi(2A_{cap}R)}\int_{B(x_0,8A_{PI}A_{cap}R)}\mathrm{d}\Gamma(f),
\end{align*}
hence
\begin{align*}
&|f(x)|^p\le2^{p-1}(|f(x)-f(y)|^p+|f(y)|^p)\\
&\le2^{p-1}\left(C_{MS}\Theta(1)\frac{\Psi(2A_{cap}R)}{\Phi(2A_{cap}R)}\int_{B(x_0,8A_{PI}A_{cap}R)}\mathrm{d}\Gamma(f)+|f(y)|^p\right).
\end{align*}
Integrating both sides over $B(x_0,A_{cap}R)$ in the variable $y$ and dividing by $V(x_0,A_{cap}R)$, we obtain that for any $x\in B(x_0,A_{cap}R)$,
$$|f(x)|^p\le 2^{p-1}\left(C_{MS}\Theta(1)\frac{\Psi(2A_{cap}R)}{\Phi(2A_{cap}R)}\int_{B(x_0,8A_{PI}A_{cap}R)}\mathrm{d}\Gamma(f)+\dashint_{B(x_0,A_{cap}R)}|f|^p\mathrm{d} m\right).$$
Therefore, $\phi\in\mathcal{F}$ is also a cutoff function for $B(x_0,R)\subseteq B(x_0,8A_{PI}A_{cap}R)$, and by the strongly local property, we have
\begin{align*}
&\int_{B(x_0,8A_{PI}A_{cap}R)}|f|^p\mathrm{d}\Gamma(\phi)=\int_{B(x_0,A_{cap}R)}|f|^p\mathrm{d}\Gamma(\phi)\\
&\le2^{p-1}\left(C_{MS}\Theta(1)\frac{\Psi(2A_{cap}R)}{\Phi(2A_{cap}R)}\int_{B(x_0,8A_{PI}A_{cap}R)}\mathrm{d}\Gamma(f)+\dashint_{B(x_0,A_{cap}R)}|f|^p\mathrm{d} m\right)\\
&\hspace{20pt}\cdot\int_{B(x_0,A_{cap}R)}\mathrm{d}\Gamma(\phi)\\
&\lesssim \left(\frac{\Psi(R)}{\Phi(R)}\int_{B(x_0,8A_{PI}A_{cap}R)}\mathrm{d}\Gamma(f)+\frac{1}{\Phi(R)}\int_{B(x_0,A_{cap}R)}|f|^p\mathrm{d} m\right)\frac{\Phi(R)}{\Psi(R)}\\
&\le\int_{B(x_0,8A_{PI}A_{cap}R)}\mathrm{d}\Gamma(f)+\frac{1}{\Psi(R)}\int_{B(x_0,8A_{PI}A_{cap}R)}|f|^p\mathrm{d} m.
\end{align*}
\end{proof}

The implication ``(1)$\Rightarrow$(2)" in Theorem \ref{thm_equiv} follows from the following two results. Assume \ref{eq_SVR}, \ref{eq_res}, then it is obvious that $\mathcal{F}\subseteq C(X)$.

\begin{lemma}\label{lem_R2cap}
Assume \ref{eq_SVR}, \ref{eq_res}. Then there exist $C>0$, $A>1$ such that for any ball $B(x_0,r)$, we have
$$R(x_0,X\backslash B(x_0,r))\ge R(B(x_0,A^{-1}r),X\backslash B(x_0,r))\ge \frac{1}{C} \frac{\Psi(r)}{\Phi(r)}.$$
Hence \ref{eq_cap} holds.
\end{lemma}

\begin{proof}
The proof is essentially the same as the proof of \cite[LEMMA 2.4]{BCK05}, \cite[Lemma 4.1]{Kum04} and \cite[Proposition 5.3]{Hu08}. Let $B=B(x_0,r)$. For any $x\in B\backslash(\frac{1}{2}B)$, there exists $h_x\in\mathcal{F}$ with $0\le h_x\le 1$ on $X$, $h_x(x_0)=1$, $h_x(x)=0$ such that
$$R(x_0,x)^{-1}\le\mathcal{E}(h_x)\le2R(x_0,x)^{-1}\le 2C_R \frac{\Phi(d(x_0,x))}{\Psi(d(x_0,x))}.$$

Let $\varepsilon\in(0,\frac{1}{4})$ be chosen later. For any $y\in B(x,\varepsilon r)$, we have
\begin{align*}
&h_x(y)=|h_x(x)-h_x(y)|\le R(x,y)^{1/p}\mathcal{E}(h_x)^{1/p}\\
&\le \left(C_R \frac{\Psi(d(x,y))}{\Phi(d(x,y))}2C_R\frac{\Phi(d(x_0,x))}{\Psi(d(x_0,x))}\right)^{1/p}\\
&\le \left(2C_R^2\Theta \left(\frac{d(x,y)}{d(x_0,x)}\right)\right)^{1/p}\le\left(2C_R^2\Theta \left(2\varepsilon\right)\right)^{1/p}.
\end{align*}
Since $\lim_{t\downarrow0}\Theta(t)=\Theta(0)=0$, there exists $\varepsilon\in(0,\frac{1}{4})$ depending only on $p, C_R, \Theta$ such that $\left(2C_R^2\Theta(2\varepsilon)\right)^{1/p}\le \frac{1}{2}$, which gives $0\le h_x\le \frac{1}{2}$ in $B(x,\varepsilon r)$. By \ref{eq_VPhi}, there exist a positive integer $N$ depending only on $C_{VR}, C_\Phi, \varepsilon$, and $x_1,\ldots,x_N\in B\backslash(\frac{1}{2}B)$ such that $B\backslash(\frac{1}{2}B)\subseteq\cup_{i=1}^NB(x_i,\varepsilon r)$. Let $h=\min_{i=1,\ldots,N}h_{x_i}$, then $0\le h\le 1$ on $X$, $h(x_0)=1$, $0\le h\le \frac{1}{2}$ in $B\backslash(\frac{1}{2}B)$. By \ref{eq_SubAdd}, we have $h\in\mathcal{F}$ and
$$\mathcal{E}(h)\le\sum_{i=1}^N\mathcal{E}(h_{x_i})\le2NC_R\Theta(1)\frac{\Phi(\frac{r}{2})}{\Psi(\frac{r}{2})}.$$
Let $g=2((h-\frac{1}{2})\vee0)$, then $g\in\mathcal{F}$ satisfies $0\le g\le 1$ on $X$, $g(x_0)=1$, $g=0$ in $B\backslash(\frac{1}{2}B)$ and $\mathcal{E}(g)\le2^p\mathcal{E}(h)$. Take $\phi\in\mathcal{F}\cap C_c(X)$ satisfying that $0\le\phi\le 1$ on $X$, $\phi=1$ on $\frac{1}{2}B$ and $\mathrm{supp}(\phi)\subseteq B$, then $\phi g\in\mathcal{F}$ satisfies $0\le\phi g\le 1$ on $X$, $(\phi g)(x_0)=1$, $\phi g=0$ on $X\backslash B$ and $\phi g=\phi g1_B=\phi g1_{\frac{1}{2}B}=g1_{\frac{1}{2}B}$, hence
\begin{align*}
&R(x_0,X\backslash B)^{-1}\le\mathcal{E}(\phi g)=\Gamma(\phi g)(X)=\Gamma(g)(\frac{1}{2}B)\le\Gamma(g)(X)=\mathcal{E}(g)\\
&\le 2^{p+1}NC_R\Theta(1)\frac{\Phi(\frac{r}{2})}{\Psi(\frac{r}{2})}\le 2^{p+1}NC_\Psi C_R\Theta(1)\frac{\Phi(r)}{\Psi(r)}=C_1\frac{\Phi(r)}{\Psi(r)},
\end{align*}
where $C_1=2^{p+1}NC_\Psi C_R\Theta(1)$, which gives $R(x_0,X\backslash B)\ge \frac{1}{C_1}\frac{\Psi(r)}{\Phi(r)}$.

Let $A>1$ be chosen later. For any $x\in A^{-1}B$, we have
\begin{align*}
&|(\phi g)(x_0)-(\phi g)(x)|\le R(x_0,x)^{1/p}\mathcal{E}(\phi g)^{1/p}\\
&\le \left(C_R \frac{\Psi(d(x_0,x))}{\Phi(d(x_0,x))}C_1 \frac{\Phi(r)}{\Psi(r)}\right)^{1/p}\le  \left(C_RC_1\Theta(A^{-1})\right)^{1/p}.
\end{align*}
Since $\lim_{t\downarrow0}\Theta(t)=\Theta(0)=0$, there exists $A>1$ depending only on $p, C_R, C_1, \Theta$ such that $\left(C_RC_1\Theta(A^{-1})\right)^{1/p}\le \frac{1}{2}$, which gives $\phi g\ge \frac{1}{2}$ in $A^{-1}B$. Let $f=2((\phi g)\wedge \frac{1}{2})$, then $f\in\mathcal{F}$ satisfies $0\le f\le 1$ on $X$, $f=1$ in $A^{-1}B$ and $f=0$ on $X\backslash B$, hence
$$\mathrm{cap}(A^{-1}B,X\backslash B)=R(A^{-1}B,X\backslash B)^{-1}\le\mathcal{E}(f)\le2^p\mathcal{E}(\phi g)\le 2^pC_1 \frac{\Phi(r)}{\Psi(r)}.$$
\end{proof}

\begin{lemma}\label{lem_R2PI}
Assume \ref{eq_SVR}, \ref{eq_res}. Then there exist $C>0$, $A>1$ such that for any $u\in\mathcal{F}$, for any $x,y\in X$ with $x\ne y$, we have
\begin{equation}\label{eq_R2PI}
|u(x)-u(y)|^p\le C \frac{\Psi(d(x,y))}{\Phi(d(x,y))}\int_{B(x,Ad(x,y))}\mathrm{d}\Gamma(u).
\end{equation}
Hence \ref{eq_PI} holds.
\end{lemma}

\begin{proof}
Firstly, we prove (\ref{eq_R2PI}). The proof is essentially the same as the proof of \cite[LEMMA 2.5]{BCK05}. Let $A>1$ be chosen later. Let $r=d(x,y)$. Without loss of generality, we may assume that $u(x)=1$, $u(y)=0$ and we only need to show that
$$\int_{B(x,Ar)}\mathrm{d}\Gamma(u)=\Gamma(u)(B(x,Ar))\gtrsim \frac{\Phi(r)}{\Psi(r)}.$$
By the Markovian property, we may assume that $0\le u\le 1$ on $X$. By Lemma \ref{lem_R2cap}, there exists $v\in\mathcal{F}$ with $0\le v\le 1$ on $X$, $v(x)=1$ and $v=0$ on $X\backslash B(x,Ar)$ such that
$$R(x,X\backslash B(x,Ar))^{-1}\le\mathcal{E}(v)\le2 R(x,X\backslash B(x,Ar))^{-1}\le 2C_1 \frac{\Phi(Ar)}{\Psi(Ar)},$$
where $C_1$ is the constant in Lemma \ref{lem_R2cap}. Let $w=u\wedge v$, then $w(x)=1$, $w(y)=0$ and $w=0$ on $X\backslash B(x,Ar)$. By \ref{eq_SubAdd}, we have $w\in\mathcal{F}$ and
\begin{align*}
&R(x,y)^{-1}\le \mathcal{E}(w)=\Gamma(w)(X)=\Gamma(w)(B(x,Ar))\le\Gamma(u)(B(x,Ar))+\Gamma(v)(B(x,Ar))\\
&\le\Gamma(u)(B(x,Ar))+\mathcal{E}(v)\le\Gamma(u)(B(x,Ar))+2C_1 \frac{\Phi(Ar)}{\Psi(Ar)}.
\end{align*}
By \ref{eq_res}, we have
$$R(x,y)^{-1}\ge \frac{1}{C_R} \frac{\Phi(r)}{\Psi(r)},$$
hence
$$\frac{1}{C_R} \frac{\Phi(r)}{\Psi(r)}\le\Gamma(u)(B(x,Ar))+2C_1 \frac{\Phi(Ar)}{\Psi(Ar)}\le\Gamma(u)(B(x,Ar))+2C_1\Theta(A^{-1})\frac{\Phi(r)}{\Psi(r)}.$$
Since $\lim_{t\downarrow0}\Theta(t)=\Theta(0)=0$, there exists $A>1$ depending only on $C_R, C_1, \Theta$ such that $2C_1\Theta(A^{-1})\le \frac{1}{2C_R}$, then $\Gamma(u)(B(x,Ar))\ge \frac{1}{2C_R}\frac{\Phi(r)}{\Psi(r)}$, which is our desired result.

Secondly, we prove \ref{eq_PI}. For any ball $B=B(x_0,R)$, for any $u\in\mathcal{F}$, by (\ref{eq_R2PI}),
\begin{align*}
&\int_B|u-u_B|^p\mathrm{d} m\\
&\le\int_B \left(\dashint_B|u(x)-u(y)|m(\mathrm{d} y)\right)^pm(\mathrm{d} x)\le \int_B\dashint_B|u(x)-u(y)|^pm(\mathrm{d} y)m(\mathrm{d} x)\\
&\le\frac{1}{m(B)}\int_B\int_B \left(C \frac{\Psi(d(x,y))}{\Phi(d(x,y))}\int_{B(x,Ad(x,y))}\mathrm{d}\Gamma(u)\right)m(\mathrm{d} y)m (\mathrm{d} x)\\
&\le C\Theta(1)\frac{\Psi(2R)}{\Phi(2R)}m(B)\int_{B(x_0,4AR)}\mathrm{d}\Gamma(u)\lesssim\Psi(R)\int_{4AB}\mathrm{d}\Gamma(u).
\end{align*}
\end{proof}

We give the proof of Theorem \ref{thm_equiv} as follows.

\begin{proof}[Proof of Theorem \ref{thm_equiv}]
``(3)$\Rightarrow$(2)": It is obvious since under \ref{eq_VD}, \ref{eq_CS} implies \ref{eq_cap}.

``(2)$\Rightarrow$(3)": It follows from Proposition \ref{prop_CS}.

``(1)$\Rightarrow$(2)": Assume \ref{eq_SVR}, \ref{eq_res}, then Lemma \ref{lem_R2cap} gives \ref{eq_cap}, and Lemma \ref{lem_R2PI} gives \ref{eq_PI}.

``(2)$\Rightarrow$(1)": By \ref{eq_SVR}, \ref{eq_PI}, applying Lemma \ref{lem_MS}, for any $u\in\mathcal{F}\subseteq C(X)$, for any $x,y\in X$ with $x\ne y$, we have
$$|u(x)-u(y)|^p\le C_{MS}\frac{\Psi(d(x,y))}{\Phi(d(x,y))}\int_{B(x,16A_{PI}d(x,y))}\mathrm{d}\Gamma(u)\le C_{MS}\frac{\Psi(d(x,y))}{\Phi(d(x,y))}\mathcal{E}(u),$$
hence
$$R(x,y)\le C_{MS}\frac{\Psi(d(x,y))}{\Phi(d(x,y))}.$$
On the other hand, by \ref{eq_cap}, noting that
$$\{x\}\subseteq B(x,\frac{1}{2A_{cap}}d(x,y)),\qquad\{y\}\subseteq X\backslash B(x,\frac{1}{2}d(x,y)),$$
we have
\begin{align*}
&R(x,y)\ge R(B(x,\frac{1}{2A_{cap}}d(x,y)),X\backslash B(x,\frac{1}{2}d(x,y)))\\
&=\mathrm{cap}(B(x,\frac{1}{2A_{cap}}d(x,y)),X\backslash B(x,\frac{1}{2}d(x,y)))^{-1}\\
&\ge \left(C_{cap}\frac{V(x,\frac{1}{2A_{cap}}d(x,y))}{\Psi(\frac{1}{2A_{cap}}d(x,y))}\right)^{-1}\asymp \frac{\Psi(d(x,y))}{\Phi(d(x,y))}.
\end{align*}
\end{proof}

\section{Some results from potential theory}\label{sec_potential}

In this section, we give some necessary results from potential theory. 

\begin{lemma}\label{lem_Leibniz}
For any $f,g\in\mathcal{F}\cap L^\infty(X;m)$, for any $A\in\mathcal{B}(X)$, we have
\begin{equation}\label{eq_Leibniz}
\Gamma(fg)(A)\le2^{p-1}\left(\int_A|\widetilde{f}|^p\mathrm{d}\Gamma(g)+\int_A|\widetilde{g}|^p\mathrm{d}\Gamma(f)\right).
\end{equation}

\end{lemma}

\begin{proof}
Firstly, we prove for $f,g\in\mathcal{F}\cap C_c(X)$. For any $\varepsilon>0$, there exist finitely many disjoint Borel sets $\{B_n:n=1,\ldots,N\}$ such that $\mathrm{supp}(f)\cup\mathrm{supp}(g)\subseteq\cup_{n=1}^NB_n$ and $|f(x)-f(y)|<\varepsilon$, $|g(x)-g(y)|<\varepsilon$ for any $x,y\in B_n$, for any $n=1,\ldots,N$. For any $n$, by \ref{eq_Alg} for $\Gamma(\cdot)(A\cap B_n)$, we have
$$\Gamma(fg)(A\cap B_n)^{1/p}\le \lVert {f}\rVert_{L^\infty(A\cap B_n;m)}\Gamma(g)(A\cap B_n)^{1/p}+\lVert {g}\rVert_{L^\infty(A\cap B_n;m)}\Gamma(f)(A\cap B_n)^{1/p},$$
hence
\begin{align*}
&\Gamma(fg)(A\cap B_n)\\
&\le2^{p-1}\left(\lVert {f}\rVert_{L^\infty(A\cap B_n;m)}^p\Gamma(g)(A\cap B_n)+\lVert {g}\rVert_{L^\infty(A\cap B_n;m)}^p\Gamma(f)(A\cap B_n)\right)\\
&\le2^{p-1}\left(\int_{A\cap B_n}(|f|+\varepsilon)^p\mathrm{d}\Gamma(g)+\int_{A\cap B_n}(|g|+\varepsilon)^p\mathrm{d}\Gamma(f)\right).
\end{align*}
Taking the summation with respect to $B_n$, we have
$$\Gamma(fg)(A)\le2^{p-1}\left(\int_{A}(|f|+\varepsilon)^p\mathrm{d}\Gamma(g)+\int_{A}(|g|+\varepsilon)^p\mathrm{d}\Gamma(f)\right).$$
Since $f,g\in\mathcal{F}\cap C_c(X)$, by the dominated convergence theorem, letting $\varepsilon\downarrow0$, we have the desired result.

Secondly, we prove for $f\in\mathcal{F}\cap L^\infty(X;m)$, $g\in\mathcal{F}\cap C_c(X)$. By \cite[Proposition 8.5]{Yan25a}, there exists $\{f_n\}\subseteq\mathcal{F}\cap C_c(X)$ which is $\mathcal{E}_1$-convergent to $f$ and convergent to $\widetilde{f}$ q.e. on $X$. By replacing $f_n$ with $(f_n\vee(-\lVert {f}\rVert_{L^\infty(X;m)}))\wedge \lVert {f}\rVert_{L^\infty(X;m)}$ and applying \cite[Corollary 3.19 (b)]{KS24a}, we may assume that $\lVert {f_n}\rVert_{L^\infty(X;m)}\le \lVert {f}\rVert_{L^\infty(X;m)}$ for any $n$. By \ref{eq_Alg}, we have $\sup_n\mathcal{E}(f_ng)<+\infty$, since $\{f_ng\}$ is $L^p(X;m)$-convergent to $fg$, by \cite[Lemma 3.17]{KS24a}, we have $\{f_ng\}$ is $\mathcal{E}_1$-weakly-convergent to $fg$. By replacing $\{f_n\}$ with a sequence of finite convex combinations of $\{f_n\}$ and applying the Mazur's lemma, we may assume that $\{f_ng\}$ is $\mathcal{E}_1$-convergent to $fg$. By the first part, we have
\begin{equation}\label{eq_Leibniz1}
\Gamma(f_ng)(A)\le2^{p-1}\left(\int_A|f_n|^p\mathrm{d}\Gamma(g)+\int_A|g|^p\mathrm{d}\Gamma(f_n)\right),
\end{equation}
where $\lim_{n\to+\infty}\Gamma(f_ng)(A)=\Gamma(fg)(A)$. Since $\{|f_n|^p\}$ converges to $|\widetilde{f}|^p$ q.e. on $X$, by \cite[Proposition 8.12]{Yan25a}, we have $\{|f_n|^p\}$ converges to $|\widetilde{f}|^p$ $\Gamma(g)$-a.e. in $X$, by the dominated convergence theorem, we have $\lim_{n\to+\infty}\int_A|f_n|^p\mathrm{d}\Gamma(g)=\int_A|\widetilde{f}|^p\mathrm{d}\Gamma(g)$. By \cite[Proposition 4.9]{KS24a}, we have $\lim_{n\to+\infty}\int_A|g|^p\mathrm{d}\Gamma(f_n)=\int_A|g|^p\mathrm{d}\Gamma(f)$. Letting $n\to+\infty$ in (\ref{eq_Leibniz1}), we have the desired result.

Finally, we prove for $f,g\in\mathcal{F}\cap L^\infty(X;m)$. Similar to the proof in the second part, there exists $\{g_n\}\subseteq\mathcal{F}\cap C_c(X)$ with $\lVert {g_n}\rVert_{L^\infty(X;m)}\le \lVert {g}\rVert_{L^\infty(X;m)}$ for any $n$, such that $\{g_n\}$ is $\mathcal{E}_1$-convergent to $g$, $\{g_n\}$ converges to $\widetilde{g}$ q.e. on $X$, and $\{fg_n\}$ is $\mathcal{E}_1$-convergent to $fg$. By the second part, we have
\begin{equation}\label{eq_Leibniz2}
\Gamma(fg_n)(A)\le2^{p-1}\left(\int_A|\widetilde{f}|^p\mathrm{d}\Gamma(g_n)+\int_A|g_n|^p\mathrm{d}\Gamma(f)\right),
\end{equation}
where $\lim_{n\to+\infty}\Gamma(fg_n)(A)=\Gamma(fg)(A)$. Since $\{|g_n|^p\}$ converges to $|\widetilde{g}|^p$ q.e. on $X$, by \cite[Proposition 8.12]{Yan25a} again, we have $\{|g_n|^p\}$ converges to $|\widetilde{g}|^p$ $\Gamma(f)$-a.e. in $X$, by the dominated convergence theorem, we have $\lim_{n\to+\infty}\int_A|g_n|^p\mathrm{d}\Gamma(f)=\int_A|\widetilde{g}|^p\mathrm{d}\Gamma(f)$. By \cite[Proposition 4.9]{KS24a} again, we have $\lim_{n\to+\infty}\int_A|\widetilde{f}|^p\mathrm{d}\Gamma(g_n)=\int_A|\widetilde{f}|^p\mathrm{d}\Gamma(g)$. Letting $n\to+\infty$ in (\ref{eq_Leibniz2}), we have the desired result.
\end{proof}

\begin{lemma}\label{lem_calFU}
For any bounded open set $U$, we have
$$\mathcal{F}(U)=\left\{u\in\mathcal{F}:\widetilde{u}=0\text{ q.e. on }X\backslash U\right\}.$$
\end{lemma}

\begin{proof}
Let $\mathcal{G}=\left\{u\in\mathcal{F}:\widetilde{u}=0\text{ q.e. on }X\backslash U\right\}$, then by \cite[Corollary 8.7]{Yan25a}, we have $\mathcal{G}$ is a closed subspace of $(\mathcal{F},\mathcal{E}^{1/p}_1)$. It is obvious that $\mathcal{F}(U)\subseteq\mathcal{G}$. We only need to show that $\mathcal{G}\subseteq\mathcal{F}(U)$. We proceed as in the proof of \cite[Lemma 5.43]{BB11}. For notational convenience, we may assume that all functions in $\mathcal{F}$ are quasi-continuous.

For any $u\in\mathcal{G}$, by \cite[Corollary 3.19 (a)]{KS24a}, we may assume that $u$ is non-negative bounded, that is, $0\le u\le M$ $m$-a.e. in $X$, where $M$ is some positive constant. For any $k\ge1$, there exists an open set $G_k$ with $\mathrm{cap}_1(G_k)<\frac{1}{k}$ such that $u|_{X\backslash G_k}$ is continuous and $u=0$ on $X\backslash(U\cup G_k)$, see \cite[Section 8]{Yan25a} for the definition of $\mathrm{cap}_1$. By \cite[Lemma 8.2]{Yan25a}, there exists $e_k\in\mathcal{F}$ with $0\le e_k\le1$ $m$-a.e. in $X$, $e_k=1$ $m$-a.e. on $G_k$ such that $\mathcal{E}_1(e_k)=\mathrm{cap}_1(G_k)$.

For any $\varepsilon>0$, we claim that $\{x\in X\backslash G_k:u(x)\ge\varepsilon\}\subseteq U$ is compact. Indeed, since $u|_{X\backslash G_k}$ is continuous, we have $(u|_{X\backslash G_k})^{-1}((-\infty,\varepsilon))=\{x\in X\backslash G_k:u(x)<\varepsilon\}$ is a relatively open subset of $X\backslash G_k$, then $G_k\cup\{x\in X\backslash G_k:u(x)<\varepsilon\}$ is an open subset of $X$, which gives $\{x\in X\backslash G_k:u(x)\ge\varepsilon\}$ is a closed subset of $X$. Since $u=0$ on $X\backslash(U\cup G_k)$, we have $\{x\in X\backslash G_k:u(x)\ge\varepsilon\}\subseteq U$. Since $U$ is bounded, we have $\{x\in X\backslash G_k:u(x)\ge\varepsilon\}$ is a bounded closed set, which is also compact due to the fact $(X,d)$ is complete locally compact.

For any $\varepsilon>0$, let $u_\varepsilon=(u-\varepsilon)\vee0$, then by \cite[Corollary 3.19 (a)]{KS24a}, we have $u_\varepsilon\in\mathcal{F}$ and $\lim_{\varepsilon\downarrow0}\mathcal{E}_1(u_\varepsilon-u)=0$.  Let $u_{\varepsilon,k}=u_\varepsilon-e_ku_\varepsilon$, then $\mathrm{supp}(u_{\varepsilon,k})\subseteq \{x\in X\backslash G_k:u(x)\ge\varepsilon\}$, by the above claim, we have $\mathrm{supp}(u_{\varepsilon,k})\subseteq U$ is compact. By \ref{eq_Alg}, we have $u_{\varepsilon,k}\in\mathcal{F}$. By Lemma \ref{lem_Leibniz}, we have
\begin{align*}
&\mathcal{E}(e_ku_\varepsilon)=\Gamma(e_ku_\varepsilon)(X)\le2^{p-1}\left(\int_X|e_k|^p\mathrm{d}\Gamma(u_\varepsilon)+\int_X|u_\varepsilon|^p\mathrm{d}\Gamma(e_k)\right)\\
&\le2^{p-1}\left(\int_X|e_k|^p\mathrm{d}\Gamma(u_\varepsilon)+M^p\mathcal{E}_1(e_k)\right)\le2^{p-1}\left(\int_X|e_k|^p\mathrm{d}\Gamma(u_\varepsilon)+M^p \frac{1}{k}\right).
\end{align*}
Since $\{e_k\}$ is $\mathcal{E}_1$-convergent to $0$, by \cite[Corollary 8.7]{Yan25a}, there exists a subsequence, still denoted by $\{e_k\}$, such that $\{e_k\}$ converges to $0$ q.e. on $X$, by \cite[Proposition 8.12]{Yan25a}, $\{e_k\}$ also converges to $0$ $\Gamma(u_\varepsilon)$-a.e. in $X$. Since $0\le e_k\le1$ q.e. on $X$, which is also $\Gamma(u_\varepsilon)$-a.e. in $X$, by the dominated convergence theorem, we have $\lim_{k\to+\infty}\int_X|e_k|^p\mathrm{d}\Gamma(u_\varepsilon)=0$. Hence $\lim_{k\to+\infty}\mathcal{E}(e_ku_\varepsilon)=0$. Moreover, since
$$\int_X|e_ku_\varepsilon|^p\mathrm{d} m\le M^p\int_X|e_k|^p\mathrm{d} m\le M^p\mathcal{E}_1(e_k)\le M^p \frac{1}{k}\to0,$$
we have $\lim_{k\to+\infty}\mathcal{E}_1(e_ku_\varepsilon)=0$, that is, $\{u_{\varepsilon,k}\}_{k}$ is $\mathcal{E}_1$-convergent to $u_\varepsilon$.

Since $u_{\varepsilon,k}\in\mathcal{F}$ satisfies $0\le u_{\varepsilon,k}\le M$ $m$-a.e. in $X$, by \cite[Corollary 3.19 (b)]{KS24a}, there exists $\{v_{\varepsilon,k,n}\}_n\subseteq\mathcal{F}\cap C_c(X)$ with $0\le v_{\varepsilon,k,n}\le M$ on $X$ for any $n$, such that $\{v_{\varepsilon,k,n}\}_n$ is $\mathcal{E}_1$-convergent to $u_{\varepsilon,k}$. Since $\mathrm{supp}(u_{\varepsilon,k})\subseteq U$ is compact, there exists $\phi_{\varepsilon,k}\in\mathcal{F}\cap C_c(X)$ satisfying $0\le\phi_{\varepsilon,k}\le1$ on $X$, $\phi_{\varepsilon,k}=1$ on $\mathrm{supp}(u_{\varepsilon,k})$ and $\mathrm{supp}(\phi_{\varepsilon,k})\subseteq U$, then $\phi_{\varepsilon,k} v_{\varepsilon,k,n}\in\mathcal{F}\cap C_c(U)$ and $u_{\varepsilon,k}-\phi_{\varepsilon,k} v_{\varepsilon,k,n}=\phi_{\varepsilon,k}(u_{\varepsilon,k}-v_{\varepsilon,k,n})$. Moreover
\begin{align*}
&\int_X|u_{\varepsilon,k}-\phi_{\varepsilon,k} v_{\varepsilon,k,n}|^p\mathrm{d} m=\int_X|\phi_{\varepsilon,k}(u_{\varepsilon,k}-v_{\varepsilon,k,n})|^p\mathrm{d} m\\
&\le\int_X|u_{\varepsilon,k}-v_{\varepsilon,k,n}|^p\mathrm{d} m\le\mathcal{E}_1(u_{\varepsilon,k}-v_{\varepsilon,k,n})\to0,
\end{align*}
and by Lemma \ref{lem_Leibniz}
\begin{align*}
&\mathcal{E}(u_{\varepsilon,k}-\phi_{\varepsilon,k} v_{\varepsilon,k,n})=\mathcal{E}(\phi_{\varepsilon,k}(u_{\varepsilon,k}- v_{\varepsilon,k,n}))\\
&\le2^{p-1}\left(\int_{X}|\phi_{\varepsilon,k}|^p\mathrm{d}\Gamma(u_{\varepsilon,k}-v_{\varepsilon,k,n})+\int_X|u_{\varepsilon,k}-v_{\varepsilon,k,n}|^p\mathrm{d}\Gamma(\phi_{\varepsilon,k})\right)\\
&\le2^{p-1}\left(\int_{X}\mathrm{d}\Gamma(u_{\varepsilon,k}-v_{\varepsilon,k,n})+\int_X|u_{\varepsilon,k}-v_{\varepsilon,k,n}|^p\mathrm{d}\Gamma(\phi_{\varepsilon,k})\right)\\
&\le2^{p-1}\left(\mathcal{E}_1(u_{\varepsilon,k}-v_{\varepsilon,k,n})+\int_X|u_{\varepsilon,k}-v_{\varepsilon,k,n}|^p\mathrm{d}\Gamma(\phi_{\varepsilon,k})\right).
\end{align*}
Since $\{v_{\varepsilon,k,n}\}_n$ is $\mathcal{E}_1$-convergent to $u_{\varepsilon,k}$, by \cite[Corollary 8.7]{Yan25a}, there exists a subsequence, still denoted by $\{v_{\varepsilon,k,n}\}_n$, such that $\{v_{\varepsilon,k,n}\}_n$ converges to $u_{\varepsilon,k}$ q.e. on $X$. By \cite[Proposition 8.12]{Yan25a}, we have $\{v_{\varepsilon,k,n}\}_n$ converges to $u_{\varepsilon,k}$ $\Gamma(\phi_{\varepsilon,k})$-a.e. in $X$. Since $0\le u_{\varepsilon,k}\le M$ q.e. on $X$, which is also $\Gamma(\phi_{\varepsilon,k})$-a.e. in $X$, and $v_{\varepsilon,k,n}\in\mathcal{F}\cap C_c(X)$ satisfies $0\le v_{\varepsilon,k,n}\le M$ on $X$, by the dominated convergence theorem, we have $\lim_{n\to+\infty}\int_X|u_{\varepsilon,k}-v_{\varepsilon,k,n}|^p\mathrm{d}\Gamma(\phi_{\varepsilon,k})=0$. Hence $\lim_{n\to+\infty}\mathcal{E}(u_{\varepsilon,k}-\phi_{\varepsilon,k} v_{\varepsilon,k,n})=0$, which gives $\{\phi_{\varepsilon,k}v_{\varepsilon,k,n}\}_n$ is $\mathcal{E}_1$-convergent to $u_{\varepsilon,k}$.

In summary, for any $\delta>0$, firstly, there exists $\varepsilon>0$ such that $\mathcal{E}_1(u_\varepsilon-u)^{1/p}<\frac{\delta}{3}$, secondly, there exists $k\ge1$ such that $\mathcal{E}_1(u_\varepsilon-u_{\varepsilon,k})^{1/p}<\frac{\delta}{3}$, finally, there exists $n$ such that $\mathcal{E}_1(u_{\varepsilon,k}-\phi_{\varepsilon,k} v_{\varepsilon,k,n})^{1/p}<\frac{\delta}{3}$, hence $\mathcal{E}_1(u-\phi_{\varepsilon,k} v_{\varepsilon,k,n})^{1/p}<\delta$, where $\phi_{\varepsilon,k} v_{\varepsilon,k,n}\in\mathcal{F}\cap C_c(U)$, that is, $u$ is the $\mathcal{E}_1$-limit of a sequence in $\mathcal{F}\cap C_c(U)$, or equivalently, $u\in\mathcal{F}(U)$.
\end{proof}

\bibliographystyle{plain}
%\bibliography{/Users/meng/Dropbox/myref}

\end{document}